\documentclass[aap,MSNbibl,seceqn,dvips]{arximspdf}

\doi{10.1214/12-AAP879} 
\volume{23}
\issue{4}
\pubyear{2013}
\firstpage{1506}
\lastpage{1543}

\makeatletter
\newcommand{\rrvert}{\vert}
\newcommand{\llvert}{\vert}
\newcommand{\eqref}[1]{(\ref{#1})}
\newtheorem{thmm}{Theorem}[section]

\newtheorem{prop}{Proposition}[section]
\newproclaim{definition}{Definition}[section]
\newtheorem{lemma}{Lemma}[section]
\newproclaim{remark}{Remark}[section]
\newproclaim{ex}{Example}[section]
\newcommand{\veps}{\varepsilon}
\newcommand{\rmd}{\mathrm{d}}
\newcommand{\rmi}{\mathrm{i}}
\newcommand{\rme}{\mathbf{e}}
\newcommand{\whH}{\widehat{H}}
\newcommand{\whV}{\widehat{V}}
\newcommand{\whk}{\widehat{\kappa}}
\newcommand{\pibar}{\overline{\Pi}}
\newcommand{\Vbar}{\overline{V}}
\newcommand{\Xbar}{\overline{X}}
\newcommand{\Ybar}{\overline{Y}}
\newcommand{\wbar}{\overline{w}}
\newcommand{\mubar}{\overline{\mu}}
\newcommand{\R}{\mathbb{R}}
\newcommand{\tux}{\tau(u-x)}
\newcommand{\tu}{\tau(u)}
\newcommand{\XD}{X\circ{\theta_{\tau(u-x)}}}
\newcommand{\X}{X_{[0,\tau(u-x))}}
\newcommand{\Xt}{X_{[0,t)}}
\newcommand{\wt}{w_{[0,t)}}
\newcommand{\tG}{\tilde{G}}
\newcommand{\PuT}{P^{(u,T)}}
\newcommand{\tZ}{\tilde{Z}}
\newcommand{\tlW}{\tilde{W}}
\newcommand{\ga}{\alpha}
\newcommand{\gl}{\lambda}
\newcommand{\gk}{\kappa}
\newcommand{\gt}{\theta}
\newcommand{\tH}{\tilde{H}}
\newcommand{\tgL}{\tilde{\Lambda}}
\newcommand{\tgF}{\tilde{\Phi}}
\newcommand{\gb}{\beta}

\makeatother

\begin{document}
\begin{frontmatter}

\title{Convolution equivalent L\'evy processes and first~passage times}
\runtitle{Convolution equivalent L\'evy processes}

\begin{aug}
\author[A]{\fnms{Philip S.} \snm{Griffin}\corref{}\ead[label=e1]{psgriffi@syr.edu}}
\runauthor{P. S. Griffin}
\affiliation{Syracuse University}
\address[A]{Department of Mathematics\\
Syracuse University\\
Syracuse, New York 13244-1150\\
USA\\
\printead{e1}} 
\end{aug}

\received{\smonth{7} \syear{2011}}
\revised{\smonth{5} \syear{2012}}

%
\begin{abstract}
We investigate the behavior of L\'{e}vy processes with convolution
equivalent L\'evy measures, up to the time of
first passage over a high level $u$. Such problems arise naturally in
the context of insurance risk where $u$ is the initial reserve.
We obtain a precise asymptotic estimate on the probability of first
passage occurring by time $T$. This result is then used to study the
process conditioned on first passage by time $T$. The existence of a
limiting process as $u\to\infty$ is demonstrated, which leads to
precise estimates for
the probability of other events
relating to first passage, such as the overshoot. A discussion of these
results, as they relate to insurance risk, is also given.
\end{abstract}

%

\end{frontmatter}

\section{Introduction}\label{s1}

Let $X=(X_{t})_{t \geq0}$,
be a L\'{e}vy process
with characteristics $(\gamma, \sigma^2, \Pi_X)$. Thus the
characteristic function of $X$
is given by the L\'{e}vy--Khintchine
representation, $Ee^{\mathrm{i}\theta X_{t}} = e^{t \Psi_X(\theta)}$,
where
\[
\label{lrep} \Psi_X(\theta) = \rmi\theta\gamma-
\sigma^2\theta^2/2+ \int_{\R}
\bigl(e^{\rmi\theta x}-1- \rmi\theta x \mathbf{1}_{\{|x|<1\}} \bigr)
\Pi_X(\rmd x) \qquad\mbox{for } \theta\in\mathbb{R}.
\]

Historically, a number of different types of L\'evy processes have
arisen in the context of stochastic modeling depending on the
phenomenon under investigation. This has motivated the detailed study
of several different classes of processes. In this paper we will
investigate one such class, those with convolution equivalent L\'evy
measure. This class has recently been proposed as a model for insurance
risk, although its study certainly predates that.
Convolution equivalent distributions were first introduced by
Chistiakov~\cite{Ch} and later by Chover, Ney and Wainger~\cite{CNW}.
Their properties have been investigated by several authors including
\cite{C,EG,kl,Pakes} and~\cite{Pakes2}
where background and further information on this class of distributions
can be found.
We will restrict ourselves to the nonlattice case,
with
the understanding that the alternative can be handled by obvious
modifications.
A distribution $F$
on $ [0, \infty)$ with tail $\overline{F}$ belongs to
the \textit{class} $\mathcal{L}^{(\alpha)}$, $\alpha\ge0$,
if $\overline{F}(u)>0$ for all $u>0$ and
\[
\lim_{u \to\infty} \frac{\overline{F}(u+x)}{\overline{F}(u)} = e^{-\alpha x} \qquad
\mbox{for } x \in(-\infty, \infty).
\]
${F}$ belongs to
the \textit{class} $\mathcal{S}^{(\alpha)}$ if, in addition,
%
\begin{equation}
\label{S2} \lim_{u \to\infty} \frac{\overline
{F^{2*}}(u)}{\overline{F}(u)} \qquad\mbox{exists and is
finite},
\end{equation}
where $F^{2*}=F * F$.
When $F\in\mathcal{S}^{(\alpha)}$,
\begin{equation}
\label{expmomS} \delta_{\alpha}^F:= \int_{[0, \infty)}
e^{\alpha x}F(\rmd x)<\infty,
\end{equation}
and the limit in \eqref{S2} is given by $ 2\delta_{\alpha}^F$.
Distributions in $\mathcal{S}^{(0)}$ are called subexponential, and
those in $\mathcal{S}^{(\alpha)}$, $\ga>0$, are called convolution
equivalent with index $\alpha>0$. The class $\mathcal{S}^{(\alpha)}$
has several nice properties, including
closure under tail equivalence, that is, if $F\in\mathcal{S}^{(\alpha
)}$ and $G$ is a distribution on $ [0, \infty)$ for which
\begin{equation}
\label{tailequiv} \lim_{u\to\infty}\frac{\overline{G}(u)}{\overline{F}(u)}=c \qquad\mbox{for some }
c\in(0,\infty),
\end{equation}
then $G\in\mathcal{S}^{(\alpha)}$. This particular property 
also holds, trivially, for the class~$\mathcal{L}^{ (\alpha)}$.

The right tail of any L\'evy measure, which is
nonzero on an interval $[x_0,\infty)$, $x_0>0$, may be taken as the
tail of a distribution function on $[x_0,\infty)$, after
\mbox{renormalization}.
With this convention, we say that the L\'evy measure (or its tail)
is in $\mathcal{S}^{ (
\alpha
)}$, respectively, $\mathcal{L}^{ (
\alpha
)}$, if this is true of the
corresponding renormalized tail. By closure under tail equivalence,
this does not depend on the choice of $x_0$. A~convolution equivalent
L\'evy process is one for which $\overline{\Pi}^+_{X} \in\mathcal
{S}^{(\alpha)}$
for some $\alpha>0$, where ${\Pi}^+_{X}$ is the restriction of ${\Pi
}_{X}$ to $(0,\infty)$ and, as above, $\overline{\Pi}^+_{X}$ denotes
its tail.
Examples include, for appropriate choices of parameters,
the CGMY, generalized inverse Gaussian (GIG) and generalized hyperbolic
(GH) processes.

Let
\begin{equation}
\label{tau} \tau(u)=\inf\{t>0\dvtx X_t>u\}
\end{equation}
denote the first passage time over level $u$.
The behavior of
\begin{equation}
\label{rat} \lim_{u\to\infty}\frac{P(\tau(u)<T)}{{\pibar}^+_{X}(u)}
\end{equation}
has been investigated under various conditions on ${\Pi}_{X}$, for
example, when ${\pibar}^+_{X}(u)$ is regularly varying\vadjust{\goodbreak} (see
Berman~\cite{berm} and Marcus~\cite{M}), and more generally when
${\pibar}^+_{X}(u)$ is subexponential (see Rosi\'nski and Samordnitsky~\cite{RS}).
In the case of interest in this paper, when $X$ is
convolution equivalent, Braverman and Samordnitsky~\cite{BS} proved
that the limit in \eqref{rat} exists but were unable to identify its
value. Later, Braverman~\cite{br} obtained a complicated description
of the limit, which unfortunately lends little practical insight as to
its actual value. Albin and Sund\'en~\cite{AS} gave a much simpler
proof of existence, but again their characterization of the limit is
highly inexplicit.
When $T=\infty$,
Kl\"{u}ppelberg, Kyprianou and Maller
\cite{kkm} were able to evaluate the limit in \eqref{rat} under the
additional assumption $ Ee^{\alpha X_1}<1$. As will become apparent in
Section~\ref{s4} (see Remark~\ref{rmkinf}), when this condition fails
the limit in \eqref{rat} is infinite for $T=\infty$.

The assumption $Ee^{\alpha X_1}<1$ was introduced in~\cite{kkm} in the
context of modeling insurance risk. Here $u$ represents the initial
reserve and $X$ the excess in claims over premium. Ruin occurs when $X$
exceeds $u$.
Our first result, which evaluates the limit in \eqref{rat}, may thus
be viewed in this context as providing a sharp asymptotic estimate for
the probability of ruin in finite time.
%
\begin{thmm}
Assume $\pibar_X^+\in\mathcal{S}^{(\alpha)}$,
then
%
\begin{equation}
\label{rt71} \lim_{u\to\infty}\frac{P(\tu<T )}{\pibar_X^+(u)} = \int_{[0,T)}e^{\psi(\ga) t}Ee^{\alpha\Xbar_{T-t}}
\,\rmd t,
\end{equation}
where $\Xbar_t=\sup_{0\le s\le t}X_s$ and $\psi(\ga) = \ln
Ee^{\alpha X_1}$.
\end{thmm}
%
The limit in \eqref{rt71} is finite, since $Ee^{\alpha\Xbar
_T}<\infty$ for every $T<\infty$ when $\pibar_X^+\in\mathcal
{S}^{(\alpha)}$ (see Lemma~\ref{EXbar}). It yields a simple and
transparent formula which 
allows further investigation of the limit as a function of $T$, as will
be illustrated in Section~\ref{s4}.
Formally, setting $\ga=0$, \eqref{rt71} reduces to the subexponential
result of~\cite{RS}. However, our interest here is in the convolution
equivalent case, so throughout this paper, it will be tacitly assumed,
without further mention, that $\ga>0$.

Building on work begun in~\cite{GM2} in the $Ee^{\alpha X_1}<1$ and
$T=\infty$ case, we investigate not only when, but how first passage
occurs in finite time, that is, what do sample paths look like that
result in first passage by time $T$?
Our main result is a functional limit theorem yielding an asymptotic
description of the process conditioned on $\tu<T$ as $u\to\infty$.
Roughly speaking, the conditioned
process behaves like an Esscher transform $Z$ of $X$ up to independent
time $\tau$ when it jumps from $Z_{\tau-}$ to a neighborhood of $u$.
Let its position after the jump be $u+W_0$. If $W_0>0$ the conditioned
process then behaves like $X$ started at $u+W_0$. If $W_0\le0$, the
conditioned process $X-u$ behaves like $X$ started at $W_0$ and
conditioned on $\tau(0)<T-\tau$. The precise descriptions of $Z, \tau
$ and $W_0$ are contained in \eqref{Z}, \eqref{tau0} and
\eqref{Wtau} and the functional limit theorem in Theorem~\ref{THM1b}.\vadjust{\goodbreak}
This result may be used to obtain precise asymptotic estimates for the
probability of many other events relating to first passage. As one
example, we
derive the joint limiting distribution of the first passage time and
the overshoot of $X$ conditional on $\tu<T$ (see Theorems~\ref{JOS}
and~\ref{JOS1}).
It will be clear from this example how other limiting distributions
relating to first passage may be found.
Previous work in this area has been restricted to the $T=\infty$ and
$Ee^{\ga X_1}<1$ case. Our results are the first that we are aware of
that considers the finite time horizon problem and removes the
condition that $Ee^{\ga X_1}<1$. The case $Ee^{\ga X_1}=1$ is of
particular interest, being the classical Cram\'er--Lundberg condition.
This is discussed further in Sections~\ref{s5} and~\ref{s7}.

We conclude the \hyperref[s1]{Introduction} with a brief outline of the paper.
Section~\ref{s2} contains various notation and introduces two measures related
to the description of the limiting process given above. Section \ref
{s3} adapts a convergence result from~\cite{GM2} in the $T=\infty$
case to the $T<\infty$ case. Section~\ref{s4} then contains the proof
of \eqref{rt71}. A~further discussion of the meaning of \eqref{rt71}
in the context of insurance risk is given in Section~\ref{s5}. Section
\ref{s6} contains the functional limit theorem and Section~\ref{s7}
applies it to the overshoot. Finally, the \hyperref[sA]{Appendix} justifies several
formulas used in the paper relating to the measures introduced in
Section~\ref{s2}.




\section{Notation}\label{s2}

We follow much of the notation laid out in~\cite{GM2}. This is briefly
summarized in the next few paragraphs for the convenience of the
reader. Let $E=\R\cup\{\Delta\}$ where $\Delta$ is a cemetery state.
Define a metric $d$ on $E$ by
\[
d(x,y)= \cases{|x-y|\wedge1, & \quad $x,y\in\mathbb R$,\vspace*{2pt}
\cr
1, &\quad $x\in
\mathbb R, y=\Delta$,\vspace*{2pt}
\cr
0, & \quad$x=y=\Delta.$ }
\]
Thus $\Delta$ is an isolated point and for $x, y\in\R$, $|x-y|\to0$
if and only if $d(x,y)\to0$.
Let $D$ be the Skorohod space of functions on $[0,\infty)$, taking
values in the metric space $E$, and which are
right-continuous with left limits.
Let
\[
\tau_z=\tau_z(w)=\inf\{t>0\dvtx w_t>z\},
\qquad \tau_\Delta= \tau_\Delta(w) = \inf\{t>0\dvtx
w_t=\Delta\}.
\]
Thus, in the notation of \eqref{tau}, $\tau(z)=\tau_z(X)$. To avoid
any possible confusion
we reserve the notation $\tau(z)$ exclusively for $\tau_z(X)$.
When considering the passage time of a process other than $X$, say $W$,
we will write $\tau_z(W)$.

For a given function
$w=(w_t)_{t\ge0}\in D$,
and
$r\ge0$, let
$w_{[0,r)}=(w_{[0,r)}(t))_{t\ge0}\in D$ denote the killed path
\[
w_{[0,r)} (t)=\cases{w_t, &\quad $0\le t<r$, \vspace*{2pt}
\cr
\Delta,&\quad $t\ge r.$}
\]
Observe that for any $t\ge0$ and $w\in D$
%
\[
\tau_\Delta(\wt)=t\qquad \mbox{if }\tau_\Delta(w)
\ge t.\vadjust{\goodbreak}
\]
For $x\in E$ let $c^x \in D$ be the constant path $c^x_t=x$ for all
${t\ge0}$. If $w, w'\in D$, then $w-w'$ denotes the path in $D$ given by
\[
\bigl(w-w'\bigr)_t= \cases{ w_t-w_t',
&\quad$\mbox{if } t<\tau_\Delta(w)\wedge \tau_\Delta
\bigl(w'\bigr)$,\vspace*{2pt}
\cr
\Delta,&\quad $\mbox{otherwise}.$ }
\]

It is convenient to assume that $X$ is given as the coordinate process
on~$D$.
The usual right-continuous completion of the filtration generated by
the coordinate maps will be denoted by $\{\mathcal{F}_t\}_{t\ge0}$.
$P_z$ denotes the probability measure induced on
$\mathcal{F}=\bigvee_{t\ge0} \mathcal{F}_t$
by the L\'evy process starting at $z\in\R$.
We usually write just $P$ for~$P_0$.
The shift operators $\theta_t\dvtx D\to D$, $t\ge0$,
are defined by $(\theta_t(w))_s=w(t+s)$.

Let $\mathcal{B}$ denote the Borel sets on $\R$ and $\mathcal
{B}([0,\infty))$ the
Borel sets on $[0,\infty)$. Let $\mathcal{D}= D\otimes[0,\infty)
\otimes(-\infty,\infty)$,
and for $T\in(0,\infty)$, set $\mathcal{D}_T= D\otimes[0,T)
\otimes(-\infty,\infty)$.
For $K\in(-\infty,\infty]$ and $x\in[0,\infty]$, define measures
$\mu_{K}$ and $\nu_x$ on $\mathcal{F}\otimes\mathcal{B}([0,\infty
))\otimes\mathcal{B}$ by
\[
\mu_K(\rmd w, \rmd t, \rmd\phi) = I(\phi< K)e^{\alpha\phi}P(\Xt
\in\rmd w, X_{{t}-}\in\rmd\phi) \,\rmd t
\]
and
\[
\nu_x \bigl(\rmd w', \rmd r, \rmd z \bigr) = I(z>-x)
\alpha e^{-\alpha z}\,\rmd z P_z \bigl(X\in\rmd w',
\tau(0) \in\rmd r \bigr).
\]
We will write $\mu$ and $\nu$ for $\mu_\infty$ and $\nu_\infty$,
respectively. The \hyperref[sA]{Appendix} contains a brief discussion of these measures,
and several formulas involving them, which will be used in the body of
the paper. Their probabilistic meaning will be discussed below after
some preliminary observations have been made.

Without any assumptions on the L\'evy process, $\mu_K$ and $\nu_x$
may be infinite measures, but on $\mathcal{D}_T$ they are finite if
$K<\infty$ and $x<\infty$, respectively. This is because
\begin{equation}
\label{mkt} \mu_K(\mathcal{D}_T)=\int
_0^T E \bigl(e^{\alpha X_t};
X_t<K \bigr) \,\rmd t
\end{equation}
and
\begin{eqnarray}
\label{nkt} \nu_x(\mathcal{D}_T)&=&\int
_{z>-x} \alpha e^{-\alpha z} P_z \bigl(\tau
(0)< T \bigr) \,\rmd z
\nonumber
\\
&=&1+\int_{0<z<x} \alpha e^{\alpha z} P \bigl(\tau(z)< T
\bigr) \,\rmd z
\\
&=&E \bigl(e^{\alpha\Xbar_T}; \Xbar_T\le x \bigr) +
e^{\alpha
x}P(\Xbar_T>x).
\nonumber
\end{eqnarray}
%
Here, and elsewhere, we make use of the fact that $X_t=X_{t-}$ a.s. for
every $t>0$.
From \eqref{mkt} and \eqref{nkt} we also see that $\mu$ and $\nu$
are finite on $\mathcal{D}_T$ whenever $Ee^{\alpha\Xbar_T}<\infty$.
This condition clearly implies $Ee^{\alpha X_1}< \infty$, and, as we
show below, is equivalent to it. This will allow us to conclude that
$\mu$ and $\nu$ are finite on $\mathcal{D}_T$ when $\pibar_X^+\in
\mathcal{S}^{(\alpha)}$.

Let $(L_t)_{t\ge0}$ be the local time of $X$ at its maximum and $H$
the corresponding ascending ladder height process\vadjust{\goodbreak} (see~\cite{bert,doneystf} or
\cite{kypbook}). The renewal function of $H$ is
\[
V(z)= \int_{t\ge0} P(H_t\le z;
t<L_\infty) \,\rmd t,
\]
with associated renewal measure $V(\rmd z)$. When $X_t\to-\infty$
a.s., $L_\infty$ has an exponential distribution with some parameter
$q>0$, $V$ is a finite measure of mass $q^{-1}$ and the following
version of the \textit{Pollacek--Khintchine formula} holds
(
see~\cite{bert}, Proposition VI.17) for $z\ge0$,
\begin{equation}
\label{PK} P \bigl(\tau(z)<\infty \bigr)=q\Vbar(z),
\end{equation}
where
\[
\Vbar(z)=\int_{x>z}V(\rmd x).
\]
Thus
\begin{equation}
\label{PK2} Ee^{\ga\Xbar_\infty}=q \int e^{\alpha z}V(\rmd z).
\end{equation}

%
\begin{lemma}\label{EXbar} If $Ee^{\alpha X_1}< \infty$ then
$Ee^{\alpha\Xbar_T}<\infty$ for every $T<\infty$. If, in addition,
$Ee^{\alpha X_1}< 1$, then $Ee^{\alpha\Xbar_\infty}<\infty$. The
condition $Ee^{\alpha X_1}< \infty$ holds when $\pibar_X^+\in
\mathcal{S}^{(\alpha)}$.
\end{lemma}

\begin{pf}
Assume $Ee^{\alpha X_1}<1$. Then $X_t\to-\infty$ a.s. since
$e^{\alpha X_t}$ is a nonnegative supermartingale. Further, by \cite
{kkm}, Proposition 5.1,
\[
\int_z e^{\alpha z}V(\rmd z)<\infty
\]
(their condition $\Pi_X^+\neq0$ is not needed for this). Hence,
$
Ee^{\ga\Xbar_\infty}<\infty
$
by \eqref{PK2}.

Now assume $1\le Ee^{\alpha X_1}<\infty$. Then we may choose $\delta
>0$ so that\break $Ee^{\alpha(X_1-\delta)}<1$. In that case $Y_t=X_t-\delta
t$ is a L\'evy process with
$Ee^{\alpha Y_1}<1$. Since $\Xbar_T\le\Ybar_T+\delta T$, it then
follows that $Ee^{\alpha\Xbar_T}<\infty$.

Finally, if $\pibar_X^+\in\mathcal{S}^{(\alpha)}$, then $\int_{[1,
\infty)} e^{\alpha x}\Pi_X(\rmd x)<\infty$ by \eqref{expmomS},
and so $Ee^{\alpha X_1}<\infty$ by~\cite{sato}, Theorem 25.17.
\end{pf}

It will be convenient to introduce measures $\mu_K^T$ and $\nu_x^T$
on $\mathcal{D}$ defined by
\[
\mu_K^T(\cdot)=\mu_K(\cdot
\cap\mathcal{D}_T ),\qquad \nu_x^T(\cdot)=
\nu_x(\cdot\cap\mathcal{D}_T ).
\]
From the above discussion, these are finite measures if $K$ and $x$ are
finite, or
if $Ee^{\alpha X_1}< \infty$.
Observe that
%
\begin{eqnarray}
\label{muTKnuTK} \mu_K^T(\rmd w, \rmd t,
\rmd\phi) &=& I(t<T)\mu_K(\rmd w, \rmd t, \rmd\phi)\nonumber\\
&=&I \bigl(
\tau_\Delta(w)<T \bigr)\mu_K(\rmd w, \rmd t, \rmd\phi),
\nonumber
\\[-8pt]
\\[-8pt]
\nonumber
\nu_x^T \bigl(\rmd w', \rmd r,
\rmd z \bigr) &=& I(r<T)\nu_x \bigl(\rmd w', \rmd r,
\rmd z \bigr)\\
&=&I \bigl(\tau_0 \bigl(w' \bigr)<T \bigr)
\nu_x \bigl(\rmd w', \rmd r, \rmd z \bigr).\nonumber
\end{eqnarray}
The first equalities are trivial and the second follow from Lemma~\ref{AL2}.
The marginal measures will be denoted in the obvious way, for example,
%
\begin{eqnarray}
\label{marg} %
\mu_K^T(\rmd w) &=& I \bigl(
\tau_\Delta(w)<T \bigr)\mu_K(\rmd w)
\nonumber
\\
&= &\int_0^T \int_{\phi< K}e^{\alpha\phi}P(
\Xt\in\rmd w, X_{{t}-}\in\rmd\phi)\,\rmd t,
\nonumber
\\[-8pt]
\\[-8pt]
\nonumber
\nu_x^T(\rmd r) &= &I(r<T)\nu_x(
\rmd r)
\\
&=& I(r<T)\int_{z>-x}\alpha e^{-\alpha z}\,\rmd z
P_z \bigl(\tau(0)\in \rmd r \bigr).
\nonumber
\end{eqnarray}
This minor abuse of notation should not cause any confusion.

The precise probabilistic meaning of $\mu$ and $\nu$ in the $\mathcal
{S}^{(\alpha)}$ case can now be given. Define processes $\tZ$ and
$\tlW$ by
\begin{equation}
\label{Ztil} P(\tZ\in\rmd w)=\frac{\mu^T(\rmd w)}{\mu^T(\mathcal{D})}= \frac
{1}{\mu^T(\mathcal{D})}\int
_{0}^{T} E\bigl(e^{\alpha X_{t-}};\Xt\in\rmd w
\bigr) \,\rmd t
\end{equation}
and
\begin{eqnarray}
\label{Wtil} P\bigl(\tlW\in\rmd w'\bigr)&=&\frac{\nu^T(\rmd
w')}{\nu^T(\mathcal{D})}
\nonumber
\\[-8pt]
\\[-8pt]
\nonumber
&=&
\frac{1}{\nu^T(\mathcal{D})}\int_{(-\infty, \infty)}\alpha e^{-\alpha z}P_z
\bigl(X\in\rmd w', \tau(0)<T\bigr) \,\rmd z.
\end{eqnarray}
It will be shown that $\tZ$ is an Esscher transform of $X$ killed at
an independent time~$\tau$ where
\[
P(\tau\in\rmd t)= \frac{\mu^T(\rmd t)}{\mu^T(\mathcal{D})}=\frac
{I(t<T)e^{\psi(\ga) t}\,\rmd t}{\mu^T(\mathcal{D})},
\]
while $\tlW$ is the process $X$ conditioned on $\tau(0)<T$, and
started with initial distribution
\[
P(\tlW_0\in\rmd z)= \frac{1}{\nu^T(\mathcal{D})}\alpha e^{-\alpha
z}P_z
\bigl(\tau(0)<T\bigr) \,\rmd z.
\]
Roughly speaking, in terms of the description of the limiting
conditioned process given in the \hyperref[s1]{Introduction,} $\mu$ describes the
behavior of the conditioned process prior to the time of the jump into
the neighborhood of $u$, and $\nu$ describes the behavior after the jump.




\section{Preliminary convergence result}\label{s3}

In this section we prove a preliminary convergence result describing
the behavior of the process for large $u$ when it jumps from a
neighborhood of the origin into a neighborhood of $u$ before time $T$,
and then passes over level $u$ before a further time $T$.
With this aim in mind, we begin by introducing a broad class of
functions to which this and other convergence results apply.\vadjust{\goodbreak}

Let $H\dvtx D\otimes D\to\R$ be measurable with respect to the product
$\sigma$-algebra and set
\[
G(w,z)= E_z \bigl[H(w,X);\tau(0) <\infty \bigr],\qquad
w\in D, z\in\R.
\]
We denote by $\mathcal{H}$ the class of such functions $H$ which satisfy
%
\begin{eqnarray}
\label{f1}&\displaystyle H \bigl(w,w' \bigr)e^{\gt w_{\tau_{\Delta}-}I(w_{\tau
_{\Delta}-}\le0)}
\qquad\mbox{is bounded for some $\gt\in[0,\ga)$};&
\\
\label{f2}& G(w,\cdot) \qquad \displaystyle\mbox{is continuous a.e. on $(-
\infty,\infty)$ for every $w\in D$}.&
\end{eqnarray}
%
For $T>0$, let $\mathcal{H}_T$ be the class of functions $H$ for which
$H(w,w')I(\tau_0(w')< T)\in\mathcal{H}$.
Conditions \eqref{f1} and \eqref{f2} hold, for example,
if $H$ is bounded and continuous in the product Skorohod topology on
$D\otimes D$. More general conditions on $H$, which
ensure that \eqref{f2} holds, are given below.
Taking $\gt>0$ in \eqref{f1} allows for certain unbounded functions
$H$.

The following result is the starting point of our investigation. It is
a consequence of
\cite{GM2}, Remark 4.1.
Let
\[
A(u,x,T)= \bigl\{\tux<T, \tu-\tux<T \bigr\}.
\]

%
\begin{thmm}\label{THM1}
Assume $\pibar_X^+\in\mathcal{L}^{(\alpha)}$ and
fix $T>0$, $x \in[0, \infty)$ and $K\in(-\infty,\infty)$.
Then, for any $H\in\mathcal{H}_T$,
%
\begin{eqnarray}
\label{thm1} && \lim_{u\to\infty}\frac{E[H(\X,\XD-c^u);X_{{{\tux}-}}< K,A(u,x,T)]}{\pibar^+_X(u)}
\nonumber
\\[-8pt]
\\[-8pt]
\nonumber
&&\qquad =\int_{{D}\otimes{D}}H \bigl(w,w'
\bigr) \mu_K^T(\rmd w)\otimes\nu_x^T
\bigl(\rmd w' \bigr).
\end{eqnarray}
\end{thmm}

\begin{pf} Fix $T>0$, $x \in[0, \infty)$ and $K\in(-\infty,\infty)$.
We first note that the limiting expression is finite, since by \eqref
{f1}, for some constant $C$ and some $\gt\in[0,\ga)$,
\begin{eqnarray}
\label{dom1} &&\int_{w\in{D}}\int
_{w'\in{D}} \bigl|H\bigl(w,w'\bigr) \bigr|\mu^T_K(
\rmd w)\nu^T_x\bigl(\rmd w'\bigr)
\nonumber
\\
&&\qquad\le C \int_{w\in{D}}\int_{w'\in{D}}
e^{-\gt w_{\tau_{\Delta
}-}I(w_{\tau_{\Delta}-}\le0)}\mu^T_K(\rmd w)\nu^T_x
\bigl(\rmd w'\bigr)
\nonumber
\\[-8pt]
\\[-8pt]
\nonumber
&&\qquad\le C \int_{w\in{D}}\int_{w'\in{D}}
\bigl(1+ e^{-\gt w_{\tau_{\Delta
}-}} \bigr)\mu^T_K(\rmd w)
\nu^T_x\bigl(\rmd w'\bigr)
\\
&&\qquad= C \nu_x(\mathcal{D}_T)\int
_0^T E \bigl(e^{\ga X_{t-}}
\bigl(1+e^{-\gt
X_{t-}}\bigr);X_{t-}<K \bigr)\,\rmd t <\infty,
\nonumber
\end{eqnarray}
where the last equality follows from \eqref{Am1}.

For $w, w'\in D, z>-x, \phi<K$ and $t\ge0$ let
\begin{eqnarray}
\label{t}\qquad  \tH\bigl(w,w'\bigr)&=&H\bigl(w,w'
\bigr)I\bigl(\tau_\Delta(w)<T\bigr)I \bigl(\tau_0
\bigl(w'\bigr)<T \bigr),
\nonumber
\\
\tG(w,z)&=& E_z \bigl[\tH(w,X);\tau(0) <\infty \bigr]\nonumber\\
&=& I \bigl(
\tau_\Delta (w)<T \bigr)E_z \bigl[H(w,X);\tau(0) <T \bigr],
\\
\qquad\tgL_u(w, \phi) &=& \int_{z>-x}\tG(w,z)
\frac{\Pi_X^+(u-\phi+\rmd
z)}{\pibar_X^+(u)},
\nonumber\\
\tgF_u(t)&=&\int_{w\in{D}}
\int_{\phi<K}\tgL_u(w,\phi)P \bigl(\Xt\in \rmd w,
X_{{t}-}\in\rmd\phi;\tux \ge t \bigr).
\nonumber
\end{eqnarray}
Then trivially $\tH\in\mathcal{H}$. Next, note that
$\tilde\mu_K$ and $\tilde\nu_x$ in~\cite{GM2}, Remark 4.1, are
simply~$\mu_K$ and $\nu_x$, respectively. Thus by \eqref{muTKnuTK},
\begin{equation}
\label{fut2} \tH\bigl(w,w'\bigr)\tilde
\mu_K(\rmd w)\tilde\nu_x\bigl(\rmd w'
\bigr)=H\bigl(w,w'\bigr) \mu_K^T(\rmd w)
\nu_x^T\bigl(\rmd w'\bigr). %
\end{equation}
In particular, by \eqref{dom1},
\begin{equation}
\label{fut1} \int_{w\in{D}}\int_{w'\in{D}}
\bigl|\tH\bigl(w,w'\bigr) \bigr|\tilde\mu_K(\rmd w)\tilde
\nu_x\bigl(\rmd w'\bigr)<\infty. %
\end{equation}
It then follows from~\cite{GM2}, Remark 4.1, that
\begin{eqnarray}
\label{fut} &&\lim_{u\to\infty}\frac{E[\tH(\X,\XD-c^u);X_{{{\tux}-}}< K, \tu
<\infty]}{\pibar_X^+(u)}
\nonumber
\\[-8pt]
\\[-8pt]
\nonumber
&&\qquad =\int
_{{D}\otimes{D}} \tH\bigl(w,w'\bigr) \tilde
\mu_K(\rmd w)\otimes\tilde\nu_x\bigl(\rmd
w'\bigr),
\end{eqnarray}
if $\tgF_u$, $u\ge u_0$, are dominated by an integrable function on
$[0,\infty)$, for some $u_0<\infty$. Once this is checked, the proof
will be complete since \eqref{fut} is the same as \eqref{thm1}
because of \eqref{fut2} and the observation that 
\[
\tau_\Delta(\X)=\tux,\qquad \tau_0 \bigl(
\XD-c^u \bigr)=\tu-\tux
\]
on $\{\tu<\infty\}$.

To prove the required domination, we modify an argument from the proof
of~\cite{GM2}, Theorem 4.1.
Fix $\veps>0$ so that $\gt+\veps\le\ga$, and write
\[
\frac{\Pi_X^+(u-\phi+\rmd z)}{\pibar_X^+(u)} =\frac{\Pi_X^+(u-\phi+\rmd z)}{\pibar_X^+ (u-\phi-x)}\frac{\pibar_X^+(u-\phi
-x)}{\pibar_X^+ (u)}.
\]
If $u>\phi+ x$, the first term in the product is a probability measure
on $(-x,\infty)$, while for the second,
by a version of Potter's bounds for regularly varying functions
(
see~\cite{BGT}, Theorem 1.5.6(ii))
there exists an $A=A_\veps$ so that
%
\begin{equation}
\label{Pibd} \frac{\pibar_X^+(u-\phi-x)}{\pibar_X^+ (u)}\le A \bigl[e^{(\alpha-\veps)(\phi+x)}\vee e^{(\alpha+\veps)(\phi+x)}
\bigr]
\end{equation}
if $u\ge1$ and $\phi+x\le u-1$.
Thus if $u_0=:(K+x+1)\vee1$, then by \eqref{f1}, \eqref{t} and~\eqref{Pibd},
for some constant $C$ depending on $H, K$, $x$,
$\ga$ and $\veps,$
%
\begin{eqnarray}
\label{Lub} \sup_{u\ge u_0}\tgL_u(w,\phi)\le C I \bigl(
\tau_\Delta(w)<T \bigr)e^{(\alpha-\veps)\phi} e^{-\gt w_{\tau
_{\Delta}-}I(w_{\tau_{\Delta}-}\le0)}
\nonumber
\\[-8pt]
\\[-8pt]
\eqntext{\mbox{all
$w\in D$, $\phi< K$},}
\end{eqnarray}
where we have used that $e^{2\veps K}e^{(\alpha-\veps)\phi}\ge
e^{(\alpha+\veps)\phi}$ if $\phi<K$ when applying \eqref{Pibd}.
In particular, for every $t\ge0$,
\begin{eqnarray*}
&&\sup_{u\ge u_0}\tgL_u(\Xt,X_{{t}-})I(X_{{t}-}<K)
\\
&&\qquad\le CI(t<T) \bigl[e^{(\ga-\veps-\gt)X_{{t}-}}I(X_{{t}-}\le
0)+e^{(\ga-\veps)X_{{t}-}}I(0<X_{{t}-}<K) \bigr]
\\
&&\qquad\le C_1I(t<T),
\end{eqnarray*}
where $C_1=C(1+ e^{(\ga-\veps) K})$, since $\ga-\veps-\gt\ge0$.
Thus for $u\ge u_0$
\[
\tgF_u(t) =E \bigl[\tgL_u(
\Xt,X_{{t}-});X_{t-}<K, \tux\ge t \bigr]\le
C_1I(t<T). %
\]
Hence, $\tgF_u$ for $u\ge u_0$ are dominated, and the proof is complete.
\end{pf}

Conditions on $H$ that ensure $H\in\mathcal{H}$ are discussed in
\cite{GM2}. In particular, by~\cite{GM2}, Proposition 4.2, if $H$
satisfies \eqref{f1}, and for all $w\in D$ and $ z\in\R$,
\begin{equation}
\label{ctsa}\quad \lim_{\veps\downarrow0}H\bigl(w,w'-c^\veps
\bigr)= H\bigl(w,w'\bigr) \qquad\mbox{a.s. $P_z\bigl(
\rmd w'\bigr)$ on $ \bigl\{\tau_0\bigl(w'
\bigr)<\infty \bigr\}$,}
\end{equation}
then $H\in\mathcal{H}$. Observe that the function $H_1(w,w')=I(\tau_0(w')<T)$ satisfies~\eqref{ctsa}, because
\[
\label{tcont1} \tau_0\bigl(w'-c^\veps
\bigr)=\tau_\veps\bigl(w'\bigr)\downarrow
\tau_0\bigl(w'\bigr) \qquad\mbox{as } \veps\downarrow0
\]
on $\{\tau_0(w')<\infty\}$ by right-continuity.
Since the class of functions satisfying \eqref{ctsa} is clearly closed
under products, it follows that if $H$ satisfies \eqref{f1} and \eqref
{ctsa}, then $H\in\mathcal{H}_T$ for every $T>0$. For example, if $H$
is bounded and $H(w,\cdot)$ is continuous in any of the usual Skorohod
topologies for every $w\in D$, then $H$ satisfies~\eqref{ctsa} and
hence, $H\in\mathcal{H}_T$ for every $T>0$. Thus $\mathcal{H}_T$ is
a broad class, containing essentially all functions that are likely to
be of interest.

In subsequent sections, we will investigate convergence of 
the first passage time and the overshoot. Similar methods could be
applied to other variables related to first passage, such as, for
example, the undershoot or the time of the maximum prior to first passage.
In applying Theorem~\ref{THM1} to the first passage time and the
overshoot, the following class of functions will prove useful. Let
$f\dvtx \R^{4}\to\R$
be a bounded Borel function which is
jointly continuous in the last
two
arguments
and set
%
\begin{equation}
\label{HcH} H\bigl(w,w'\bigr)=f \bigl(
\tau_\Delta(w), \wbar_{\tau_\Delta-}, \tau_0
\bigl(w'\bigr), w'_{\tau_0} \bigr) 
\end{equation}
%
on $\{\tau_\Delta(w)<\infty, \tau_0(w')<\infty\}$, where $\wbar_t=\sup_{0\le s\le t}w_s$. Since we only consider such $H$ on this
set, it's definition elsewhere does not much matter.\vadjust{\goodbreak} For completeness
though, here and below, we take any $H$ of the form \eqref{HcH} to be
$0$ off this set.
Then by~\cite{GM2}, Proposition 5.1, $H$ satisfies \eqref{ctsa} and
hence, $H\in\mathcal{H}_T$ for every $T>0$.






\section{First passage time}\label{s4}

To study the first passage time, we begin by applying Theorem \ref
{THM1} with
\[
H\bigl(w,w'\bigr)=h\bigl(\tau_\Delta(w),
\tau_0\bigl(w'\bigr)\bigr),
\]
where $h\dvtx \R^2\to\R$ is a bounded Borel function such that $h(t,\cdot
)$ is continuous for every $t\ge0$. Then $H$ is of the form \eqref
{HcH}, and so $H\in\mathcal{H}_T$. Thus if $\pibar_X^+\in\mathcal
{L}^{(\alpha)}$, then by \eqref{thm1}
%
\begin{eqnarray}
\label{rt1} &&\lim_{u\to\infty}\frac{E[h(\tux, \tu
-\tux),
X_{{{\tux}-}} < K,A(u,x,T)]}{\pibar_X^+(u)}
\nonumber\\
&&\qquad = \int_{{D}\otimes{D}}h \bigl(\tau_\Delta(w),
\tau_0 \bigl(w' \bigr) \bigr) \mu_K^T(
\rmd w)\otimes\nu_x^T \bigl(\rmd w' \bigr)
\\
&&\qquad = \int_0^\infty\int_0^\infty
h(t,r)\mu_K^T( \rmd t)\nu_x^T(
\rmd r),\nonumber
\end{eqnarray}
where the final equality comes from applying \eqref{ltb12a} to the
positive and negative parts of $h$.
In particular,
\begin{eqnarray*}
\label{rt1a} &&\lim_{u\to\infty}\frac{P(\tux\in
\rmd t, \tu
-\tux\in\rmd r, X_{{{\tux}-}} < K)}{\pibar_X^+(u)}
\\
&&\qquad = \mu_K^T(\rmd t)\otimes \nu_x^T(
\rmd r)
\end{eqnarray*}
in the sense of weak convergence of measures on $[0,T)\otimes[0,T)$.

To obtain the asymptotic behavior of $P(\tau(u)<T)$ when $\pibar_X^+\in\mathcal{S}^{(\alpha)}$ we need a version of \eqref{rt1} in
which $X_{{{\tux}-}} < K$ is replaced by $\Xbar_{{{\tux}-}} < K$.
For this we introduce a measure $\mubar_K^T$ on $[0,\infty)$ by
\[
\mubar_K^T(\rmd t)= I(t<T)E
\bigl(e^{\alpha X_{t-}}; \Xbar_{t-}<K\bigr) \,\rmd t.
\]
This may be compared with the marginal measure [{cf}. \eqref{marg}]
\[
\mu_K^T(\rmd t) = I(t<T)
\mu_K(\rmd t) = I(t<T) E \bigl(e^{\alpha X_{t-}}; X_{t-}<K
\bigr)\,\rmd t.
\]

%
\begin{prop}\label{Xblt} Let $h\dvtx \R^2\to\R$ be a bounded Borel
function such that $h(t,\cdot)$ is continuous for every $t\ge0$. If
$\pibar_X^+\in\mathcal{L}^{(\alpha)}$, then
%
\begin{eqnarray}
\label{ltb}
&& \lim_{u\to\infty}\frac{E[h(\tux,\tu
-\tux
);\Xbar_{{{\tux}-}}<K,A(u,x,T)]}{\pibar_X^+(u)}
\nonumber
\\[-8pt]
\\[-8pt]
\nonumber
&&\qquad =\int_0^\infty\int_0^\infty
h(t,r)\mubar_K^T( \rmd t)\nu_x^T(
\rmd r).
\end{eqnarray}
\end{prop}

\begin{pf} Let $H(w,w')=h(\tau_\Delta(w), \tau_0(w'))I(\wbar_{\tau
_\Delta-}<K)$. Then $H$ is of the form \eqref{HcH}, and so by Theorem
\ref{THM1}, the limit in \eqref{ltb} is given by
\begin{eqnarray*}
&&\int_{{D}\otimes{D}} h \bigl(\tau_\Delta(w),
\tau_0 \bigl(w' \bigr) \bigr) I(\wbar_{\tau_\Delta-}<K)
\mu_K^T(\rmd w)\otimes\nu_x^T
\bigl(\rmd w' \bigr)
\\[-2pt]
&&\qquad =\int_0^\infty\int_0^\infty
h(t,r)\mubar_K^T( \rmd t)\nu_x^T(
\rmd r) 
\end{eqnarray*}
by \eqref{ltb11}.
\end{pf}

Recall that $\mu=\mu_\infty$ and $\nu=\nu_\infty$. 
If $Ee^{\ga X_1}<\infty$ then, as noted following~\eqref{nkt}, the
marginal measures
\begin{equation}
\label{mninf} \mu( \rmd t)=Ee^{\alpha X_t} \,\rmd t, \qquad\nu(\rmd r) = \int
_{z}\alpha e^{-\alpha z}\,\rmd z P_z\bigl(
\tau(0)\in\rmd r \bigr)
\end{equation}
are finite on $[0,T)$. Equivalently,
$\mu^T(\rmd t)=\mu(\rmd t\cap[0,T))$ and $\nu^T(\rmd r)=\nu(\rmd
r\cap[0,T))$
are finite measures.
For notational convenience, we will write $\nu(t)$ for $\nu([0,t))$
and similarly for other measures.

In the next two propositions and elsewhere, we consider limits as
\mbox{$K,x\to\infty$}. By this we will always mean that the manner in which
$K$ and $x$ approach infinity is irrelevant. In particular, they can do
so in either order.
Letting ${K,x\to\infty}$ in Propositions~\ref{rtt} and~\ref{rttt}
indicates that the most probable paths along which $X$ can reach level
$u$ by time $T$ are those in which the process jumps from a
neighborhood of $0$ to a neighborhood of $u$. This will be elucidated
upon further in Theorem
\ref{THM1b}.

%
\begin{prop}\label{rtt} If $\pibar_X^+\in\mathcal{L}^{(\alpha)}$
and $Ee^{\ga X_1}<\infty$, then
%
\begin{equation}
\label{rt21} \lim_{K,x\to\infty}\lim_{u\to\infty}\frac{P(\tu
<T, \Xbar_{{{\tux}-}} < K)}{\pibar_X^+(u)} = \int
_{[0,T)}\nu(T-t)\mu(\rmd t).
\end{equation}
\end{prop}

\begin{pf}
Fix $\veps>0$ and let $g$ be continuous with $g\equiv1$ on $[0, T]$,
$g\equiv0$ on $[T+\veps, \infty)$ and $g$ linear
on $[T,T+\veps]$. Since $\{\tu<T\}\subset A(u,x,T)$, we have by
\eqref{ltb} with
$h(t,r)=g(t+r)$, for every $x$ and $K$,
%
\begin{eqnarray}
\label{rt3} &&\limsup_{u\to\infty}
\frac{P(\tu<T,
\Xbar_{{{\tux}-}} < K)}{\pibar_X^+(u)}
\nonumber\\[-2pt]
&&\qquad \le\int_0^\infty\int_0^\infty
g(t+r)\mubar_K^T( \rmd t)\nu_x^T(
\rmd r)
\nonumber\\[-2pt]
&&\qquad \le\int_{[0,\infty)}\nu_x^T(T+ \veps-t)
\mubar_K^T(\rmd t)
\\[-2pt]
&&\qquad \le\int_{[0,T)}E \bigl(e^{\alpha X_{t-}};
\Xbar_{t-}<K \bigr)\,\rmd t\nonumber\\
&&\quad\qquad{}\times\int_{z}I(z>-x)\alpha
e^{-\alpha z}P_z \bigl(\tau(0)<T+\veps-t \bigr)\,\rmd z.\nonumber
\end{eqnarray}
Similarly,
%
\begin{eqnarray}
\label{rt31} &&\liminf_{u\to\infty}
\frac{P(\tu<T,
\Xbar_{{{\tux}-}} < K)}{\pibar_X^+(u)}
\nonumber
\\
&&\qquad \ge\int_{[0,T)}E \bigl(e^{\alpha X_{t-}};
\Xbar_{t-}<K \bigr)\,\rmd t\\
&&\qquad\quad{}\times \int_{z}I(z>-x)\alpha
e^{-\alpha z}P_z \bigl(\tau(0)<T-\veps-t \bigr)\,\rmd z.\nonumber
\end{eqnarray}
The result now follows by letting ${K,x\to\infty}$ and then
$\veps\downarrow0$ in \eqref{rt3} and \eqref{rt31} and
noting that
\[
\int_{[0,T)}\nu \bigl([0,T-t )\bigr)\mu(\rmd t)=
\int_{[0,T)}\nu \bigl([0,T-t] \bigr)\mu(\rmd t)
\]
since the integrands agree except on an at most countable set, and $\mu
$ has no atoms.
\end{pf}

%
\begin{remark} If $\pibar_X^+\in\mathcal{L}^{(\alpha)}$ but
$Ee^{\ga X_1}=\infty$, it follows from \eqref{rt31} that Proposition
$\ref{rtt}$ remains valid provided we interpret the integral in
\eqref{rt21} as infinite.
\end{remark}

If 
$Ee^{\alpha X_1}<1$ then $\mu$ and $\nu$ given by \eqref{mninf} are
finite measures on
$[0,\infty)$, since
\begin{equation}
\label{ml} \mu(\infty)=\int_0^\infty
Ee^{\alpha X_t} \,\rmd t =\int_0^\infty
\bigl(Ee^{\alpha X_1}\bigr)^t \,\rmd t <\infty
\end{equation}
and
\begin{equation}
\label{nl} \nu(\infty) = \int_{z}\alpha e^{-\alpha z}\,
\rmd z P_z(\tau_0<\infty )=Ee^{\alpha\Xbar_\infty}<\infty
\end{equation}
by \eqref{nkt} and Lemma~\ref{EXbar}.

%
\begin{prop}\label{rttt}
If $ \pibar_X^+\in\mathcal{S}^{(\alpha)}$, then
%
\begin{equation}
\label{rt51} \lim_{K,x\to\infty}\lim_{u\to\infty}\frac{P(\tu
<T, \Xbar_{{{\tux}-}} \ge K)}{\pibar_X^+(u)} = 0.
\end{equation}
\end{prop}

\begin{pf}By Lemma~\ref{EXbar}, $Ee^{\alpha X_1}<\infty$. First
assume that $Ee^{\alpha X_1}<1$. Then,
by~\eqref{rt21},
\begin{eqnarray*}\label{rt22}
\lim_{K,x\to\infty}\lim_{u\to\infty
}
\frac
{P(\tu<\infty, \Xbar_{{{\tux}-}} < K)}{\pibar_X^+(u)} &\ge& \lim_{T\to\infty} \int_{[0,T)}
\nu(T-t)\mu(\rmd t)
\\
&=&\mu(\infty)\nu(\infty).
\end{eqnarray*}
On the other hand, by Theorem 4.1 and~\cite{kkm}, Proposition 5.3,
%
\begin{equation}
\label{kkmt} \lim_{u\to\infty}\frac{P(\tu<\infty)}{\pibar_X^+(u)} = \mu(\infty)\nu(\infty),
\end{equation}
from which \eqref{rt51} immediately follows (with $T$ even replaced by
$\infty$).

Now assume $Ee^{\alpha X_1}\ge1$.
Choose $\delta>0$ so that $Ee^{\alpha(X_1-\delta)}<1$ and set
$Y_t=X_t-\delta t$. Then $\Pi_Y=\Pi_X$ and $Ee^{\alpha Y_1}<1$. Hence,
\eqref{rt51} holds with $X$ replaced by~$Y$. Next, recalling $\tau_{u}(Y)=\inf\{t>0\dvtx Y_t>u\}$, it is clear that
\begin{equation}
\label{compXY} \bigl\{\tu<T\bigr\}\subset\bigl\{\tau_{(u-x)}(Y)<T\bigr\} \qquad\mbox{for
} x\ge\delta T.
\end{equation}
%
Hence, if $x\wedge K>\delta T$,
then
\begin{eqnarray}
\label{comp} &&\bigl\{\tu<T, \Ybar_{\tau_{(u-x)}(Y)-}<K-\delta
T\bigr\}
\nonumber\\
&&\qquad = \bigl\{\tu<T, \Ybar_{\tau_{(u-x)}(Y)-}<K-\delta T, \tau_{(u-x)}(Y)<T
\bigr\}
\nonumber
\\[-8pt]
\\[-8pt]
\nonumber
& &\qquad\subset \bigl\{\tu<T, \Ybar_{\tau_{(u-x)}(Y)-}+\delta{\tau_{(u-x)}(Y)}<K
\bigr\}
\\
&&\qquad \subset\bigl\{\tu<T, \Xbar_{\tau_{(u-x)}(Y)-}<K\bigr\}.\nonumber
\end{eqnarray}
On the other hand, we trivially have
\[
X_{\tau_{(u-x)}(Y)}>Y_{\tau_{(u-x)}(Y)}\ge u-x.
\]
Thus, if additionally $u-x>K$, then $\tau_{(u-x)}(Y)=\tau(u-x)$ on $\{
\tu<T, \Ybar_{\tau_{(u-x)}(Y)-}<K-\delta T\}$. Consequently, by
\eqref{comp},
\begin{equation}\qquad
\label{XY} \bigl\{\tu<T, \Ybar_{\tau_{(u-x)}(Y)-}<K-\delta T\bigr\}\subset\bigl\{\tu<T,
\Xbar_{\tau(u-x)-}<K\bigr\}.
\end{equation}
%
Hence, from \eqref{compXY} and \eqref{XY} we may conclude that if
$x\wedge K>\delta T$ and $x+K<u$, then
\[
\bigl\{\tu<T, \Xbar_{{{\tux}-}} \ge K\bigr\}\subset\bigl\{\tau_{(u-\delta
T)}(Y)<T,
\Ybar_{\tau_{(u-x)}(Y)-} \ge K-\delta T\bigr\}.
\]
Thus by \eqref{rt51}, with $X$ replaced by $Y$, $u$ by $u-\delta T$,
$x$ by $x-\delta T$ and $K$ by $K-\delta T$, we have
\begin{eqnarray*}\label{rt8} &&\lim_{K,x\to\infty}\lim_{u\to\infty
}
\frac
{P(\tu<T, \Xbar_{{{\tux}-}} \ge K)}{\pibar_X^+(u)}
\\
&&\qquad \le e^{\ga\delta T}\lim_{K,x\to\infty}\lim_{u\to\infty}
\frac{P(\tau_{(u-\delta T)}(Y)<T, \Ybar_{\tau_{((u-\delta
T)-(x-\delta T))}(Y)-} \ge K-\delta T)}{\pibar_X^+(u-\delta T)}
\\
&&\qquad = 0,
\end{eqnarray*}
which completes the proof.
\end{pf}

%
\begin{remark} As noted in the proof, if $ \pibar_X^+\in\mathcal
{S}^{(\alpha)}$ and $Ee^{\alpha X_1}<1$ then
\[
\lim_{K,x\to\infty}\lim_{u\to\infty}\frac{P(\tu<\infty, \Xbar_{{{\tux}-}} \ge K)}{\pibar_X^+(u)} = 0.
\]
%
\end{remark}

Recall that
$\psi$, introduced in \eqref{rt71}, denotes the exponent of the mgf
of $X_1$, that is,
\[
e^{\psi(\gb)}=Ee^{\gb X_1}.
\]
Note that $\psi(\ga)<\infty$ and $\psi(\gb)=\infty$ for $\gb>
\ga$ when $ \pibar_X^+\in\mathcal{S}^{(\alpha)}$. Combining the
two previous propositions yields the following result.

%
\begin{thmm}[(Probability of ruin in finite time)]\label{THMall} Assume that
$\pibar_X^+\in\mathcal{S}^{(\alpha)}$, then
%
\begin{eqnarray}
\label{rt7} \lim_{u\to\infty}\frac{P(\tu<T
)}{\pibar_X^+(u)} &=&
\int_{[0,T)} \nu(T-t)\mu(\rmd t)
\nonumber
\\[-8pt]
\\[-8pt]
\nonumber
&=& \int_{[0,T)}e^{\psi(\ga) t}Ee^{\alpha\Xbar_{T-t}} \,\rmd t.
\end{eqnarray}
\end{thmm}

\begin{pf} The first equality follows from Propositions~\ref{rtt} and
\ref{rttt}, and the second from \eqref{mninf} and \eqref{nkt}.
\end{pf}

%
\begin{remark}\label{rmkinf}
If $\pibar_X^+\in\mathcal{S}^{(\alpha)}$ and $Ee^{\alpha X_1}\ge
1$, then trivially
\[
\lim_{T\to\infty}\int_{[0,T)}e^{\psi(\ga) t}Ee^{\alpha\Xbar_{T-t}}
\,\rmd t =\infty.
\]
Thus ${\pibar_X^+(u)}=o({P(\tu<\infty)})$ in contrast to the case
$Ee^{\alpha X_1}< 1$, where \eqref{kkmt} holds and the limit is finite
by \eqref{ml} and \eqref{nl}.
\end{remark}

The existence of the limit in \eqref{rt7} was first proved by
Braverman and Samordnitsky~\cite{BS}. Later, Braverman~\cite{br} obtained
a complicated description of the limit (see also Albin and Sund\'en
\cite{AS}, (6.1) and (6.6) and Hao and Tang~\cite{ht}, (4.8)).
Albin and Sund\'en's approach involved showing
\[
\lim_{u\to\infty}\frac{P(\Xbar_T>u )}{P(X_T>u)}= L(T)
\]
exists. Their description of 
$L(T)$ is highly inexplicit, but they were able to show $L(T)>1$ for
all $T>0$ when $X$ is not a subordinator. Since
\begin{equation}
\label{tailPi} \lim_{u\to\infty}\frac{P(X_T>u)}{\pibar_X^+(u)}=Te^{\psi(\ga) T}
\end{equation}
%
when $\pibar_X^+\in\mathcal{S}^{(\alpha)}$ (see~\cite{Sb})
it follows from \eqref{rt7} that
\[
L(T)= \frac1T \int_{[0,T)}e^{-\psi(\ga) t}Ee^{\alpha\Xbar_{t}} \,
\rmd t,
\]
providing an alternative proof that $L(T)>1$ precisely when $X$ is not
a subordinator.

We turn now to the limit in \eqref{rt7} and investigate its behavior
as a function of~$T$.
When $\psi(\ga)\ge0$, a change variable gives
\begin{equation}
\label{bpos} \int_{[0,T)}e^{\psi(\ga) t}Ee^{\alpha\Xbar_{T-t}}
\,\rmd t =e^{\psi(\ga) T}\int_{[0,T)}e^{-\psi(\ga) t}Ee^{\alpha\Xbar
_{t}}
\,\rmd t.
\end{equation}
The integral on the RHS diverges as $T\to\infty$ since $Ee^{\alpha
\Xbar_{t}}\ge Ee^{\alpha X_{t}}=e^{\psi(\ga) t}.$
To determine the correct exponential rate of growth, we note that $\ln
Ee^{\alpha\Xbar_{t}}$ is subadditive,
hence, by Fekete's lemma (\cite{HP}, Theorem 7.6.1),
\[
\lim_{t\to\infty}\frac{\ln Ee^{\alpha\Xbar_{t}}}{t}=C
\]
for some $C$, where
clearly $C\in[\psi(\ga),\infty)$.
It then easily follows that
\begin{equation}
\label{ga>} \lim_{T\to\infty}\frac1T\ln\biggl(\lim_{u\to\infty}
\frac{P(\tu<T
)}{\pibar_X^+(u)} \biggr)= C.
\end{equation}
In general, when $ \pibar_X^+\in\mathcal{S}^{(\alpha)}$, we only
know that $Ee^{\alpha X_{1}}<\infty$, but now assume for the remainder
of this paragraph that, in addition to $\psi(\ga)\ge0$, we have
\begin{equation}
\label{Cr1} E\bigl(X_1e^{\alpha X_{1}}\bigr)<\infty.
\end{equation}
This assumption arises in connection with the Cram\'er--Lundberg large
deviation estimate in the $\psi(\ga)=0$ case [see \eqref{ruinprob}
below]. Then, using Doob's $L^1$-maximal inequality (see~\cite{Dur},
Exercise 5.4.6), it is easy to check that for some constant $C<\infty$,
\begin{equation}
\label{Doob} Ee^{\alpha\Xbar_{t}}\le C(1+t)e^{\psi(\ga)t}.
\end{equation}
Hence,
\[
\lim_{t\to\infty}\frac{\ln Ee^{\alpha\Xbar_{t}}}{t}=\psi(\ga).
\]
Thus, in particular, we are able to identify the correct rate of
exponential growth as
\begin{equation}
\label{growth} \lim_{T\to\infty}\frac1T\ln\biggl(\lim_{u\to\infty}
\frac{P(\tu<T
)}{\pibar_X^+(u)} \biggr)= \psi(\ga).
\end{equation}
When $\psi(\ga)=0$, more precise information is available from \eqref
{Doob}. In the special case that $X$ is a subordinator, it follows from\vadjust{\goodbreak}
\eqref{rt7} or \eqref{tailPi} that \eqref{growth} can be sharpened to
\[
\lim_{u\to\infty}\frac{P(\tu<T )}{\pibar_X^+(u)}=Te^{\psi(\ga) T}.
\]
%

When $\psi(\ga)>0$, \eqref{bpos} may be rewritten
\begin{equation}
\label{B<1} \int_{[0,T)}e^{\psi(\ga) t}Ee^{\alpha\Xbar_{T-t}}
\,\rmd t =\psi(\ga)^{-1}e^{\psi(\ga) T}E\bigl(e^{\alpha\Xbar_{\mathbf{e}}}; \rme<
T\bigr),
\end{equation}
where $\mathbf{e}$ is exponentially distributed with parameter $\psi
(\ga)$ and independent of $X$.
For fixed $T$ this provides a formula which appears well suited to
Monte Carlo simulation. 
It gives the more precise, than \eqref{ga>}, asymptotic estimate
\[
\lim_{u\to\infty}\frac{P(\tu<T )}{e^{\psi(\ga)T}\pibar_X^+(u)}\sim
\frac{E(e^{\alpha\Xbar_{\rme}}; \rme< T)}{\psi(\ga)}\qquad \mbox{as $T\to\infty$,}
\]
where note
\[
Ee^{\alpha\Xbar_{\rme}}\ge Ee^{\alpha X_{\rme}}=\infty.
\]

When $\psi(\ga)<0$, that is, $Ee^{\ga X_1}<1$,
\[
\int_{[0,T)}e^{\psi(\ga) t}Ee^{\alpha\Xbar_{T-t}} \,\rmd t=-\psi (
\ga)^{-1}E\bigl[e^{\alpha\Xbar_{T-\rme}};\rme< T\bigr],
\]
where $\rme$ is exponentially distributed with parameter $-\psi(\ga
)$ and independent of $X$. In this case
\begin{equation}
\label{tauTu} \lim_{T\to\infty}\lim_{u\to\infty}\frac{P(\tu<T )}{\pibar_X^+(u)}=
\frac{Ee^{\alpha\Xbar_{\infty}}}{-\psi(\ga)}<\infty
\end{equation}
by Lemma~\ref{EXbar}. By way of comparison, observe that from \eqref
{ml} and \eqref{nl}, \eqref{kkmt} may be written
\begin{equation}
\label{tauuT} \lim_{u\to\infty}\lim_{T\to\infty}\frac{P(\tu<T )}{\pibar_X^+(u)}=
\frac{Ee^{\alpha\Xbar_{\infty}}}{-\psi(\ga)}.
\end{equation}

The asymptotic behavior as $T\to0$, irrespective of the value of $\psi
(\ga)$, also follows easily from \eqref{rt7}:
\[
\lim_{T\to0}\frac1T\lim_{u\to\infty}\frac{P(\tu<T )}{\pibar_X^+(u)}= 1.
\]




\section{An application to insurance risk}\label{s5}

A popular model in insurance risk is the Cram\'er--Lundberg model in which
\begin{equation}
\label{CL} X_t=\sum_1^{N_t}
U_i - pt,\vadjust{\goodbreak}
\end{equation}
where $N_t$ is a rate $\lambda$ Poisson process, and $U_i> 0$ form an
independent i.i.d. sequence. Here $p$ represents the rate of premium
inflow and $U_i$ the size of the $i$th claim. Thus $X$, called the
claim surplus process, represents the excess in claims over premium.
The insurance company starts with a positive reserve $u$, and ruin
occurs if this level is exceeded by $X$. It is assumed that $EU_1=\mu$
is finite and that $p=(1+\gt)\lambda\mu$ where $\gt>0$ is called
the safety loading. This ensures $X_t\to-\infty$, and so the
probability of eventual ruin $P(\tu<\infty)\to0 $ as $u\to\infty$.

A common assumption on $X$ is the Cram\'{e}r--Lundberg condition
\begin{equation}
\label{Cr} E e^{\nu X_1}=1 \qquad\mbox{for some } \nu> 0.
\end{equation}
The $\nu$ which satisfies \eqref{Cr} is often referred to as the
Lundberg exponent or adjustment coefficient. This condition results in
a well-known large deviation estimate for the probability of eventual ruin:
\begin{equation}
\label{ruinprob} \lim_{u\to\infty}e^{\nu u}P\bigl(\tu<\infty\bigr)=C=:
\frac
{-EX_1}{E(X_1e^{\nu X_{1}})},
\end{equation}
where $C>0$ if \eqref{Cr1} holds with $\ga$ replaced by $\nu$.

The problem of a sharp estimate for the probability of ruin in finite
time when an exponential moment exists is much more difficult. In the
special case that the claim size distribution is itself exponential, an
exact formula for $P(\tu<T)$ is available (see~\cite{Arp},
Proposition V.1.3). Other than this, little is known, although
several approximations have been proposed (see~\cite{Arp}, Chapter V).
One typical such approximation is 
the classical Segerdahl approximation; if \eqref{Cr} holds and
$E(X_1^2e^{\nu X_{1}})<\infty$, then
\begin{equation}
\label{Sap} P\bigl(\tau(u)<T\bigr)= Ce^{-\nu u}\Phi\biggl(
\frac{T-au}{b\sqrt{u}} \biggr)+o\bigl(e^{-\nu u}\bigr)
\end{equation}
uniformly in $T$, where $a$ and $b$ are known constants and $\Phi$ is
the standard normal distribution function.
Considerable care must be taken in using \eqref{Sap}.
The only time \eqref{Sap} is guaranteed to provide a valid estimate is
when $T\ge au+O(\sqrt{u})$. For~$T$ of smaller order, the estimate is
of smaller order than $e^{-\nu u}$. For example, for fixed $T$, \eqref
{Sap} gives
\begin{equation}
\label{Sap1} P\bigl(\tau(u)<T\bigr)= \frac{Cb}{a(2\pi u)^{1/2}}e^{-(\nu+{a^2}/{(2b^2)})
u+({a}/{b^2})T}+o
\bigl(e^{-\nu u}\bigr),
\end{equation}
and it is quite likely that the error term will exceed the estimate
itself. While some improvements to this estimate are possible, and
alternative approximations such as the (corrected) diffusion
approximation have been proposed, none can lay claim to giving a sharp
estimate for the probability of ruin in finite time.

Recently a more general L\'evy risk insurance model has been proposed
in which~\eqref{CL} is replaced by a spectrally positive L\'evy
process $X$, that is, $\Pi_X=\Pi_X^+$, for which $X_t\to-\infty$.
Theorem~\ref{THMall} then solves the problem of a sharp estimate for
the probability of ruin in finite time in this more general model (even
without the spectrally positive assumption)
when $\pibar_X^+\in\mathcal{S}^{(\alpha)}$:
%
\begin{equation}
\label{CEApp} P\bigl(\tau(u)<T\bigr)= \pibar_X^+(u)B(T) +o\bigl(
\pibar_X^+(u)\bigr),
\end{equation}
where
\begin{equation}
\label{B} B(T)= \int_{[0,T)}e^{\psi(\ga)({T-t}) }Ee^{\alpha\Xbar_t}
\,\rmd t.
\end{equation}
Since $B$ is continuous and the limit of monotone functions, the
estimate is uniform on compacts in $T$, and when $B$ is bounded, that
is, when $E e^{\ga X_1}<1$, the estimate is uniform over all $T$.
There seems little hope of evaluating $B$ explicitly
and so in practice numerical techniques will be needed.
One possibility is to approximate~$B$ using Monte Carlo simulation, for
which formulations like \eqref{B<1} appear well suited.
An alternative is to approximate $B$ by (numerically) inverting its
Laplace transform. For $\delta>\psi(\ga)\vee0$,
this is given by
\[
\int_{s\ge0}
e^{-\delta s}B(s) \,\rmd s=\frac1{\delta-\psi(\ga )}\int_{t\ge0}
e^{-\delta t}Ee^{\alpha\Xbar_t} \,\rmd t =\frac{Ee^{\alpha\Xbar_{ \mathbf{e}(\delta)}}}{\delta(\delta-\psi(\ga))},
\]
where $\bf e(\delta)$ is an independent exponential with parameter
$\delta$.
In the spectrally positive case this may be written equivalently as
\[
\int_{s\ge0}
e^{-\delta s}B(s) \,\rmd s =\frac{(\phi(\delta)-\ga)}{(\delta-\psi(\ga))^2\phi(\delta)},
\]
where $\phi$ is the inverse of the restriction of $\psi$ to $(-\infty
,0]$ (see~\cite{doneystf}, (4.3.7) and~(9.2.9) (which note applies to spectrally negative
processes)). It would be interesting to
investigate how successfully
these, and possibly other approximation methods, could be implemented
in concrete classes of examples, such as those mentioned in the
\hyperref[s1]{Introduction} or the GTSC class of models introduced by Hubalek and
Kyprianou~\cite{HK} and further investigated in~\cite{GMCr}.

As an illustration we compare estimates \eqref{Sap1} and \eqref{CEApp}
in the context of the Cram\'er--Lundberg model \eqref{CL} when $\pibar_X^+\in\mathcal{S}^{(\alpha)}$.
This is equivalent to the assumption
$U_1\in\mathcal{S}^{(\alpha)}$.
One may regard \eqref{CL} as giving a family of models indexed by the
premium rate. Let $p_0=\gl\mu$ be the premium rate corresponding to
zero safety loading,\vspace*{-1pt} and write $X^{(p)}_t=Y_t - pt$,
where $Y_t=\sum_1^{N_t} U_i$. Since $E e^{\ga Y_1}>1$, there is a
unique $p=p_L$ such that $E e^{\ga X^{(p_L)}_1}=1.$ Observe that
$X^{(p_L)}_t\to-\infty$ a.s. since $e^{\ga X^{(p_L)}_t}$ is a
nonnegative martingale. Thus $p_L>p_0$.
In comparing \eqref{Sap1} and \eqref{CEApp}, we consider three
different regimes for $p$. The first is large premiums; $p>p_L$. In
that case $E e^{\ga X^{(p)}_1}<1$, and since $E e^{\gb
X^{(p)}_1}=\infty$ for every $\gb>\ga$ when $\pibar_X^+\in\mathcal
{S}^{(\alpha)}$, the Lundberg exponent does not exist. Thus the
classical Segerdahl approximation has nothing to say in this case.
The second regime is when $p=p_L$. Then the Segerdahl approximation
yields \eqref{Sap1} with $\nu=\ga$. However, not surprisingly, the
estimate is of completely the wrong order, since $e^{-\gb u}=o(\pibar_X^+(u))$ for every $\gb>\ga$. Finally, the third regime is small
premiums; $p_0<p<p_L$.
In this case the Lundberg exponent $\nu$ exists and $\nu<\ga$,
however, again one can show that the estimate is of the wrong order. In
this third regime
an alternative approximation, the (corrected) diffusion approximation,
is often suggested. This is a heavy traffic limit, that is, an
approximation as $p\downarrow p_0$. Asmussen and Albrecher~\cite{Arp},
Chapter~V.6, report that numerical evidence indicates it provides quite
good estimates when $p$ is close to $p_0$, but again it cannot expect
to match the sharp estimate \eqref{CEApp} which is valid for every $p$.

In concluding this section it should be pointed out that one would not expect
the estimates for the probability of ruin in finite time that have been
proposed in the literature to be as good as \eqref{CEApp} when
$\pibar_X^+\in\mathcal{S}^{(\alpha)}$. After all, \eqref{CEApp} is
a sharp estimate derived from the additional structure resulting from
the assumption $\pibar_X^+\in\mathcal{S}^{(\alpha)}$. Given how
little is known about these ruin probabilities in general, \eqref
{CEApp} might be useful as a benchmark against which to compare these
more general approximations. We should also mention that in the
subexponential case the situation is much better understood. Then Rosi\'
nski and Samordnitsky~\cite{RS} show
\[
\lim_{u\to\infty}\frac{P(\tau(u)<T)}{{\pibar}^+_{X}(u)}=\lim_{u\to\infty}\frac{P(X_T>u)}{{\pibar}^+_{X}(u)}=T.
\]
The first equality is because ruin by time $T$ is essentially the
result of one extremely large claim which greatly exceeds $u$.
Consequently, $X$ will not have returned to level $u$ by time $T$ on
the event $\tau(u)<T$. The second equality is a direct consequence of
subexponentiality.





\section{Functional limit theorem}\label{s6}

We now address the question of how first passage occurs by time $T$, by
proving a functional limit theorem for the process conditioned on $\tau
(u)<T$ as $u\to\infty$. We begin by revisiting Theorem~\ref{THM1},
in the $\mathcal{S}^{(\alpha)}$ case, with the aid of Proposition
\ref{rttt}. This allows us to
set $K=\infty$ and take the limit as $x\to\infty$ in \eqref{thm1}.

%
\begin{thmm}\label{THM1a} Assume
$\pibar_X^+\in\mathcal{S}^{(\alpha)}$ and $H\in\mathcal{H}_T$. Then
%
\begin{eqnarray}
\label{1a} &&\lim_{x\to\infty}\lim_{u\to\infty
}
\frac{E[H(\X
,\XD-c^u);A(u,x,T)]}{\pibar^+_X(u)}
\nonumber
\\[-8pt]
\\[-8pt]
\nonumber
&&\qquad =\int_{{D}\otimes{D}} H \bigl(w,w' \bigr)
\mu^T(\rmd w)\otimes\nu^T \bigl(\rmd w'
\bigr).
\end{eqnarray}
\end{thmm}

\begin{pf}
The limit is finite since, by the same argument as in \eqref{dom1} but
with $K=x=\infty$, we obtain
%
\begin{eqnarray}
\label{dom}&& \int_{w\in{D}}\int
_{w'\in{D}} \bigl|H\bigl(w,w'\bigr)\bigr|\mu^T(
\rmd w)\nu^T\bigl(\rmd w'\bigr)
\nonumber
\\[-8pt]
\\[-8pt]
\nonumber
&&\qquad\le C \nu(\mathcal{D}_T)\int_0^T
E\bigl[e^{\ga X_{t-}}\bigl(1+e^{-\gt
X_{t-}}\bigr)\bigr]\,\rmd t <\infty,
\end{eqnarray}
where finiteness follows from Lemma~\ref{EXbar}.


Next, by \eqref{f1}, if $K\ge0$, then $|H|$ is bounded on $\{
(w,w')\dvtx w_{\tau_{\Delta}-}\ge K\}$ by some constant $C$ say. Since
$\tu<2T$ on $A(u,x,T)$, it then follows from \eqref{rt51} with $T$
replaced by $2T$, that
\begin{eqnarray*}\label{1aaa}
 &&\lim_{K, x\to\infty}
\lim_{u\to\infty
}\frac
{E[H(\X,\XD-c^u);X_{\tux-}\ge K, A(u,x,T)]}{\pibar^+_X(u)}
\\
&&\qquad = 0.
\end{eqnarray*}
Thus by \eqref{thm1}, the limit in \eqref{1a} is given by
\begin{eqnarray*}
&&\lim_{K,x\to\infty}\int_{{D}\otimes{D}} H \bigl(w,w'
\bigr) \mu^T_K(\rmd w)\otimes\nu^T_x
\bigl(\rmd w' \bigr)
\\
&&\qquad =\lim_{K,x\to\infty}\int_{\mathcal{D}\otimes\mathcal{D}} H \bigl(w,w'
\bigr)I(\phi< K)I(z>-x) \mu^T(\rmd w,\rmd t,\rmd\phi)\\
&&\qquad\hspace*{60pt}\quad{}\otimes
\nu^T \bigl(\rmd w', \rmd r, \rmd z \bigr).
\end{eqnarray*}
Now the integrand is trivially dominated by $|H|$ and
\begin{eqnarray*}
&&\int_{\mathcal{D}\otimes\mathcal{D}} \bigl|H
\bigl(w,w'\bigr)\bigr| \mu^T(\rmd w,\rmd t,\rmd\phi)\otimes
\nu^T\bigl(\rmd w', \rmd r, \rmd z\bigr)
\\
&&\qquad= \int_{w\in{D}}\int_{w'\in{D}}\bigl|H
\bigl(w,w'\bigr)\bigr|\mu^T(\rmd w)\nu^T\bigl(
\rmd w'\bigr)<\infty
\end{eqnarray*}
by \eqref{dom}. Thus by dominated convergence
\begin{eqnarray*}
&&\lim_{K,x\to\infty}\int_{{D}\otimes{D}} H \bigl(w,w'
\bigr) \mu_K^T(\rmd w)\otimes\nu_x^T
\bigl(\rmd w' \bigr)
\\
&&\qquad =\int_{\mathcal{D}\otimes\mathcal{D}} H \bigl(w,w' \bigr)
\mu^T(\rmd w,\rmd t,\rmd\phi)\otimes\nu^T \bigl(\rmd
w', \rmd r, \rmd z \bigr)
\\
&&\qquad=\int_{{D}\otimes{D}} H \bigl(w,w' \bigr)
\mu^T(\rmd w)\otimes\nu^T \bigl(\rmd w'
\bigr).
\end{eqnarray*}
\upqed\end{pf}

To give a clearer understanding of the limit in \eqref{1a}, introduce
independent \mbox{$D$-valued} random variables $\tZ$ and $\tlW$
with distributions given by \eqref{Ztil} and \eqref{Wtil}, respectively.\vadjust{\goodbreak}
Clearly $\tlW$ is the process $X$ conditioned on $\tau(0)<T$, and
started with initial distribution
%
\begin{equation}
\label{init} P(\tlW_0\in\rmd z)= \frac{1}{\nu^T(\mathcal
{D})}\alpha
e^{-\alpha z}P_z \bigl(\tau(0)<T \bigr) \,\rmd z.
\end{equation}
To give a more transparent description of $\tZ$, we first introduce
the Esscher transform $Z$ of $X$, defined as follows. Let $\mathcal
{B}([0,s])$ denote the Borel sets in $\mathbb R^{[0,s]}$. Then for any
$s\ge0$ and any $B_s\in\mathcal{B}([0,s])$,
\begin{eqnarray}\label{Z} P\bigl(\{Z_v\dvtx 0\le v\le s
\}\in B_s\bigr)=e^{-\psi(\ga) s}E\bigl(e^{\alpha X_{s}};\{
X_v\dvtx 0\le v\le s\}\in B_s\bigr).
\end{eqnarray}
%
%
Next, recalling \eqref{mninf}, let $\tau$ be independent of $Z$ with
distribution
\begin{equation}
\label{tau0} P(\tau\in\rmd t)= \frac{\mu^T(\rmd t)}{\mu^T(\mathcal{D})}=\frac
{I(t<T)e^{\psi(\ga) t}\,\rmd t}{\mu^T(\mathcal{D})}.
\end{equation}
Thus
\begin{equation}
\label{tau1} P(\tau\in\rmd t)= \frac{\psi(\ga) e^{\psi(\ga) t}\,\rmd t}{e^{\psi(\ga) T}-1}, \qquad 0\le t < T, \mbox{ if } \psi(
\ga) \neq0
\end{equation}
and
\begin{equation}
\label{tau2} P(\tau\in\rmd t)= \frac{\rmd t}{T},\qquad 0\le t < T, \mbox{ if } \psi(
\ga) = 0.
\end{equation}

%
\begin{prop}
Assume $Ee^{\ga X_1}<\infty$.
Then with $Z$ and $\tau$ as above,
\[
\bigl\{\tZ_t\dvtx t<\tau_{\Delta}(\tZ)\bigr\}
\stackrel{d} {=}\{Z_t\dvtx t<\tau\}.
\]
\end{prop}

\begin{pf} For any $B_s\in\mathcal{B}([0,s])$
\begin{eqnarray*}\label{YX}
&& P\bigl(\{\tZ_v\dvtx 0\le v
\le s\}\in B_s, s<\tau_{\Delta}(\tZ)\bigr)
\\
&&\qquad =\frac1{\mu^T(\mathcal{D})}\int_{s<t<T}E
\bigl(e^{\alpha
X_{t-}}\dvtx \{X_v\dvtx 0\le v\le s\}\in
B_s\bigr) \,\rmd t
\\
&&\qquad =\frac{E(e^{\alpha X_{s}}\dvtx \{X_v\dvtx 0\le v\le s\}\in
B_s)}{\mu^T(\mathcal{D})}\int_{s<t<T}e^{\psi(\ga) (t-s)} \,\rmd t
\\
&&\qquad =P\bigl(
\{Z_v\dvtx 0\le v\le s\}\in B_s\bigr)
\frac{1}{\mu^T(\mathcal
D)}\int_{s<t<T} e^{\psi(\ga) t} \,\rmd t
\\
& &\qquad=P\bigl(\{Z_v\dvtx 0\le v\le s\}\in B_s\bigr)P(s<
\tau)
\\
&&\qquad =P\bigl(\{Z_v\dvtx 0\le v\le s\}\in B_s, s<\tau
\bigr).
\end{eqnarray*}
\upqed\end{pf}

Thus $\tZ$ is seen to be the Esscher transform of $X$ killed at an
independent time~$\tau$ with distribution given by \eqref{tau0}.

With the previous analysis at hand, it is a relatively easy matter to
study the process $X$ conditioned on $\tu<T$, when
$\pibar_X^+\in\mathcal{S}^{(\alpha)}$. To do so, first introduce
the probability measure
\[
P^{(u,T)}( \cdot)=P\bigl( \cdot|\tu<T\bigr)\vadjust{\goodbreak}
\]
and let $E^{(u,T)}$ denote expectation with respect to $P^{(u,T)}$.
Let $Z$ and $\tau$ be distributed as above and let $(W,\tau)$ be
independent of $Z$ with joint distribution
%
\begin{eqnarray}
\label{Wtau} && P \bigl(W\in\rmd w', \tau
\in\rmd t \bigr)
\nonumber\\
&&\qquad=\frac{\mu^T(\mathcal{D})}{B(T)} \int_{(-\infty, \infty)}\alpha e^{-\alpha z}P_z
\bigl(X\in\rmd w', \tau(0)<T-t \bigr)\,\rmd z P(\tau\in\rmd t)
\\
&&\qquad=\frac{\mu^T(\mathcal{D})}{B(T)} \nu^{T-t} \bigl(\rmd w' \bigr)P(\tau
\in\rmd t),\nonumber
\end{eqnarray}
where recall $B(T)$ is given by \eqref{B}.
Observe this is a true probability distribution since, from \eqref{nkt},
\begin{eqnarray*} &&\int_{[0,T)}\int
_{D}\int_{(-\infty, \infty)}\alpha
e^{-\alpha
z}P_z\bigl(X\in\rmd w, \tau(0)<T-t\bigr)\,\rmd z P(
\tau\in\rmd t)
\\
&&\qquad=\frac1{\mu^T(\mathcal{D})}\int_{[0,T)}
Ee^{\alpha\Xbar_{T-t}} \mu(\rmd t)
\\
&&\qquad=\frac1{\mu^T(\mathcal{D})}\int_{[0,T)}
e^{\psi(\ga)
t}Ee^{\alpha\Xbar_{T-t}} \,\rmd t = \frac{B(T)}{\mu^T(\mathcal{D})}.
\end{eqnarray*}
Thus $W$ is the process $X$ conditioned on $\tau(0)<T-\tau$, and
started with initial distribution
\[
P(W_0\in\rmd z)= \frac{\mu^T(\mathcal{D})}{B(T)}\alpha
e^{-\alpha z}P_z \bigl(\tau(0)<T-\tau \bigr) \,\rmd z.
\]
In particular,
\[
P(W_0>0)= \frac{\mu^T(\mathcal{D})}{B(T)}.
\]

%
\begin{lemma} The joint distribution of $(Z_{[0, \tau)} ,W)$ is given by
\begin{equation}
\label{jtdist} P\bigl(Z_{[0, \tau)}\in\rmd w, W\in\rmd w'
\bigr)=\frac1{B(T)}\int_0^\infty
\mu^T(\rmd w,\rmd t)\nu^{T-t}\bigl(\rmd w'
\bigr).
\end{equation}
\end{lemma}

\begin{pf} First observe that by \eqref{Z}, for $0\le t< T$,
\begin{equation}
\label{Z1} P(Z_{[0,t)}\in\rmd w)=e^{-\psi(\ga)t}E\bigl(e^{\ga X_t};X_{[0,t)}
\in\rmd w\bigr).
\end{equation}
Thus by \eqref{tau0}, \eqref{Wtau}, \eqref{Z1} and independence
\begin{eqnarray*}
&& P\bigl(Z_{[0, \tau)}\in\rmd w, W\in\rmd
w'\bigr)
\\
&&\qquad=\int_t P\bigl(Z_{[0,t)}\in\rmd w, \tau\in\rmd t,
W\in\rmd w'\bigr)
\\
&&\qquad=\int_t P(Z_{[0,t)}\in\rmd w) P\bigl(\tau\in\rmd
t, W\in\rmd w'\bigr)
\\
&&\qquad=\frac{\mu^T(\mathcal{D})}{B(T)}\int_t e^{-\psi(\ga)t}E
\bigl(e^{\ga
X_t};X_{[0,t)}\in\rmd w\bigr)\nu^{T-t}\bigl(
\rmd w'\bigr)P(\tau\in\rmd t)
\\
&&\qquad=\frac1{B(T)}\int_{[0,T)} E\bigl(e^{\ga X_t};X_{[0,t)}
\in\rmd w\bigr)\nu^{T-t}\bigl(\rmd w'\bigr)\,\rmd t
\\
&&\qquad=\frac1{B(T)}\int_0^\infty\mu^T(
\rmd w,\rmd t)\nu^{T-t}\bigl(\rmd w'\bigr).
\end{eqnarray*}
\upqed\end{pf}

We are now ready to state the main result of this section.

%
\begin{thmm}[(Functional limit theorem)]\label{THM1b} Assume that
$\pibar_X^+\in\mathcal{S}^{(\alpha)}$ and
$H\in\mathcal{H}_t$ for every $t\le T$. Then
%
\begin{equation}
\qquad\lim_{x\to\infty}\lim_{u\to\infty
}E^{(u,T)}H
\bigl( \X,\XD-c^u \bigr)=EH(Z_{[0, \tau)} ,W). \label{1aa}
\end{equation}
\end{thmm}

\begin{pf} Set
\[
\tilde{H}\bigl(w,w'\bigr)=H\bigl(w,w'\bigr)I\bigl(
\tau_\Delta(w)+\tau_0\bigl(w'\bigr)<T\bigr),
\]
and fix $w_0\in D$. Then
\[
E_z\bigl[\tH(w_0,X)I\bigl(\tau(0) <T\bigr);\tau(0) <
\infty\bigr]=E_z\bigl[H(w_0,X);\tau(0) <T-
\tau_{\Delta}(w_0)\bigr].
\]
Now for every $w\in D$, and in particular for $w=w_0$, $E_z[H(w,X);\tau
(0) <T-\tau_{\Delta}(w_0)]$ is trivially continuous in $z$ if $\tau_{\Delta}(w_0)\ge T$, while if $\tau_{\Delta}(w_0)<T$, it is
continuous a.e. in $z\in(-\infty,\infty)$ since $H\in\mathcal
{H}_{T-\tau_{\Delta}(w_0)}$. Thus $\tH\in\mathcal{H}_T$ and so
\begin{eqnarray*}
&& E^{(u,T)}H\bigl(\X,\XD-c^u\bigr)
\\
&&\qquad=\frac{E\tilde{H}(\X,\XD-c^u)}{P(\tu<T)}
\\
&&\qquad=\frac{E[\tilde{H}(\X,\XD-c^u);A(u,x,T)]}{P(\tu<T)}
\\
&&\qquad\to\frac{1}{B(T)}\int_{{D}\otimes{D}} \tH\bigl(w,w'
\bigr) \mu^T(\rmd w)\otimes\nu^T\bigl(\rmd
w'\bigr) \end{eqnarray*}
as $u\to\infty$, then $x\to\infty$ by \eqref{1a} and \eqref{rt7}.
This last integral is absolutely convergent by \eqref{dom}. Hence,
applying \eqref{AlH} to the positive and negative parts of $H$, this
final expression may be rewritten as
\begin{eqnarray*}&& \frac{1}{B(T)}\int_{{D}\otimes{D}}  H
\bigl(w,w'\bigr)I\bigl(\tau_{\Delta}(w)+ \tau_0
\bigl(w'\bigr)<T\bigr) \mu^T(\rmd w)\otimes
\nu^T\bigl(\rmd w'\bigr)
\\
&&\qquad=\frac{1}{B(T)}\int_{{D}\otimes{D}} H\bigl(w,w'
\bigr)\int_{[0,\infty)} \mu^T(\rmd w,\rmd t)
\nu^{T-t}\bigl( \rmd w'\bigr)
\\
&&\qquad=EH(Z_{[0, \tau)} ,W)
\end{eqnarray*}
by \eqref{jtdist}.
\end{pf}

Thus, under $P^{(u,T)}$, for large $u$, the process $X$ can be
approximated as follows:
\begin{itemize}\label{ltproc}
\item run $Z$ for times $0\le t<\tau$;
\item run $u+W$ from time $\tau$ on, that is, at time $\tau+t$, the
value of the process is $u+W_t$.
\end{itemize}
Thus the process behaves like $Z$ up to an independent time $\tau$
when it jumps from a neighborhood of $0$ to a neighborhood of $u$. Its
position prior to the jump is $Z_{\tau-}$ and its position after is
$u+W_0$. If $W_0>0$ the process $X-u$ behaves like $X$ started at
$W_0$. If $W_0\le0$, the process $X-u$ behaves like $X$ started at
$W_0$ and conditioned on $\tau(0)<T-\tau$.




\setcounter{equation}{0}
\section{Conditional distribution of the first passage time and
overshoot}\label{s7}

As an illustration of Theorem~\ref{THM1b}, we derive the joint
limiting distribution of the first passage time $\tu$ and the
overshoot $O_u=X_{\tu}-u$, conditional on $\tu<T$.
We give two descriptions of the limit, the first in terms of the
limiting variables in Theorem~\ref{THM1b}, and the second in terms of
fluctuation quantities. This latter description allows us to relate the
limiting distribution of the overshoot when $T<\infty$ to the limiting
distribution when $T=\infty$ (as found in~\cite{kkm}) for the
$Ee^{\ga X_1}<1$ case.

%
\begin{thmm}\label{JOS} If $\pibar_X^+\in\mathcal{S}^{(\alpha)}$,
then for any bounded continuous function $g\dvtx [0,\infty)^2\to\R$,
\[
\lim_{u\to\infty}E^{(u,T)}g\bigl(O_u, \tu\bigr) =Eg
\bigl(W_{\tau_0(W)}, \tau+\tau_0(W)\bigr).
\]
Furthermore, the limiting distribution is given by
\begin{eqnarray}\label{ltpto}
&&P\bigl(W_{\tau_0(W)}\in\rmd\gamma,
\tau+\tau_0(W)\in\rmd t\bigr)
\nonumber\\
&&\qquad=\frac{1}{B(T)} \biggl(\alpha e^{-\alpha\gamma}\,\rmd\gamma\mu(\rmd t)\\
&&\qquad\hspace*{30pt}\quad{}+
 \int_{z\ge0}\alpha e^{\alpha z}\,\rmd z\int
_{0\le r\le t}\mu(\rmd t-r)P_0\bigl(X_{\tau(z)}-z
\in\rmd\gamma, \tau(z)\in\rmd r\bigr) \biggr)\nonumber
\end{eqnarray}
for $\gamma\ge0$ and $0\le t< T$.
\end{thmm}

\begin{pf}
Set
\[
H\bigl(w,w'\bigr)=g\bigl(w'_{\tau_0},
\tau_{\Delta}(w)+\tau_0\bigl(w'\bigr)\bigr).\vadjust{\goodbreak}
\]
Then for every $0\le x<u$,
\[
H\bigl(\X,\XD-c^u\bigr)=g\bigl(O_u, \tu\bigr)
\]
if $\tu<\infty$.
Further, $H\in\mathcal{H}_t$ for every $t>0$, since $H$ is of the
form \eqref{HcH}. Thus by \eqref{1aa} we obtain
\[
\lim_{u\to\infty}E^{(u,T)}g\bigl(O_u, \tu\bigr) =Eg
\bigl(W_{\tau_0(W)}, \tau+\tau_0(W)\bigr).
\]
Now by \eqref{tau0} and \eqref{Wtau},
\begin{eqnarray*}
&&P\bigl(W_{\tau_0(W)}\in\rmd\gamma, \tau\in\rmd
s, \tau_0(W)\in\rmd r\bigr)
\\
&&\qquad= \frac{I(r<T-s)}{B(T)} \int_{z}\alpha e^{-\alpha z}P_z
\bigl(X_{\tau(0)}\in \rmd\gamma, \tau(0)\in\rmd r\bigr)\,\rmd z
\mu^T(\rmd s).
\end{eqnarray*}
Hence, letting $\delta_{\{b\}}$ denote a point mass concentrated at
$b$, the limiting distribution of $(O_u, \tu)$ is given by
\begin{eqnarray*}
&&P\bigl(W_{\tau_0(W)}\in\rmd\gamma, \tau+
\tau_0(W)\in\rmd t\bigr)
\\
&&\qquad= \frac{1}{B(T)} \int_{z}\alpha e^{-\alpha z}\,
\rmd z\int_{0\le r\le
t}\mu(\rmd t-r)P_z
\bigl(X_{\tau(0)}\in\rmd\gamma, \tau(0)\in\rmd r\bigr)
\\
&&\qquad=\frac{1}{B(T)} \biggl(\int_{z>0}\alpha
e^{-\alpha z}\,\rmd z\int_{0\le
r\le t}\mu(\rmd t-r)
\delta_{\{(z,0)\}}(\rmd\gamma, \rmd r)\\
&&\hspace*{42pt}\qquad{}+
\int_{z\le0}\alpha e^{-\alpha z}\,\rmd z\int
_{0\le r\le t}\mu (\rmd t-r)P_z\bigl(X_{\tau(0)}
\in\rmd\gamma, \tau(0)\in\rmd r\bigr) \biggr)
\\
&&\qquad=\frac{1}{B(T)} \biggl(\alpha e^{-\alpha\gamma}\,\rmd\gamma\mu(\rmd t)\\
&&\hspace*{42pt}\qquad{}+
\int_{z\ge0}\alpha e^{\alpha z}\,\rmd z\int
_{0\le r\le t}\mu (\rmd t-r)P_0\bigl(X_{\tau(z)}-z
\in\rmd\gamma, \tau(z)\in\rmd r\bigr) \biggr)
\end{eqnarray*}
for $\gamma\ge0$ and $0\le t< T$.
\end{pf}

The expression for the limiting distribution in \eqref{ltpto} may also
be written in terms of fluctuation quantities by using the quintuple
law of Doney and Kyprianou~\cite{DK}.
In order to do so, we first need to introduce some further notation
which is standard in the area;
cf.~\cite{bert,doneystf,kypbook}.
Recall from Section~\ref{s2} that $(L_t)_{t\ge0}$ is the local time
of~$X$ at its maximum. Let $(L^{-1}_t,H_t)_{t \geq0}$ be the bivariate
ascending ladder process and $\Pi_{{ L}^{-1},{ H}}(\cdot,\cdot)$ its
L\'evy measure.
The bivariate renewal function of $(L^{-1},H)$ is
\begin{equation}
\label{BivV} V(s,z)= \int_{t\ge0} P\bigl(L_t^{-1}
\le s, H_t\le z; t<L_\infty\bigr) \,\rmd t,
\end{equation}
with associated renewal measure $V(\rmd s,\rmd z)$.
The L\'evy measure of $H$ will be denoted $\Pi_{{ H}}$, and its
Laplace exponent by $\gk$ where
\begin{eqnarray}
\label{kappa} \gk(\gb)&:=&-\ln\bigl(E\bigl(e^{-\gb H_1};1<L_\infty
\bigr)\bigr)
\nonumber
\\[-9pt]
\\[-9pt]
\nonumber
&\hspace*{1pt}=&q+\gb\rmd_{ H}+\int_{y\ge0}
\bigl(1-e^{-\gb y}\bigr)\Pi_{{ H}}(\rmd y).
\end{eqnarray}
Here $\rmd_{ H}\ge0$ is the drift and $q\ge0$ is the killing rate of
$H$ (see, e.g.,~\cite{kypbook}, (6.15) and (6.16)).

Define measures $\eta_\alpha^V(\rmd s)$ and $\eta_\alpha^{\Pi_{{
L}^{-1}, { H}}}(\rmd s)$ on $[0,\infty)$ by
\[
\eta_\alpha^V(\rmd s)=\int_{z}e^{\alpha z}V(
\rmd s,\rmd z),\qquad \eta_\alpha^{\Pi_{{ L}^{-1}, { H}}}(\rmd s)=\int
_{z}\bigl(e^{\alpha
z}-1\bigr)\Pi_{{ L}^{-1},{ H}}(\rmd
s,\rmd z).
\]
If $Ee^{\ga X_1}<\infty$ then, as we now show, $\eta_\alpha^{\Pi_{{
L}^{-1}, H}}$ is a finite measure, while $\eta_\alpha^V$ is finite on
compact sets.

%
\begin{prop}\label{Vig} If $Ee^{\ga X_1}<\infty$ then $\int_{z\ge1}
e^{\ga z}\Pi_H(\rmd z)<\infty$.
\end{prop}

\begin{pf}
By Vigon's \'equation amicale invers\'ee~\cite{Vig}, for $z>0$
\[
\Pi_H(\rmd z)=k\int_{y\ge0} \whV(
\rmd y)\Pi_X(y+\rmd z),
\]
where $\whV$ is the renewal function of the descending ladder height
process, and $k>0$ is a constant depending on the normalizations of the
local times.
Thus
\begin{eqnarray*}
\int_{z\ge1}
e^{\ga z}\Pi_H(\rmd z)&=&k\int_{z\ge1}
e^{\ga z} \int_{y\ge0} \whV(\rmd y)
\Pi_X(y+\rmd z)
\\[-2pt]
&=& k\int_{y\ge0} \whV(\rmd y)\int_{x\ge y+1}
e^{\ga(x-y)}\Pi_X(\rmd x)
\\[-2pt]
&\le& k\int_{y\ge0} e^{-\ga y} \whV(\rmd y)\int
_{x\ge1} e^{\ga
x}\Pi_X(\rmd x).
\end{eqnarray*}
The first integral is finite since for some $c>0$, $\whV([0,y])\le cy$
for large $y$ by~\cite{bert}, Proposition III.1, and the second is finite
by~\cite{sato}, Theorem 25.17. 
\end{pf}

Finiteness of $\eta_\alpha^{\Pi_{{ L}^{-1}, H}}$ follows immediately
from Proposition~\ref{Vig}. For $\eta_\alpha^V$, Proposition \ref
{Vig} and~\cite{sato}, Theorem 25.17, imply that when $Ee^{\ga
X_1}<\infty$, we also have $E(e^{\ga H_1};1<L_\infty)<\infty$.
Hence, by dominated convergence,
for $a$ sufficiently large, $E(e^{-aL^{-1}_1+\ga H_1};1<L_\infty)<1$.
Thus, for such $a$,
\begin{eqnarray*} \int_se^{-a s}
\eta_\alpha^V(\rmd s) &=&\int_s\int
_ze^{-as+\ga z}\int_{t\ge0} P
\bigl(L_t^{-1}\in\rmd s,H_t\in \rmd z,
t<L_\infty\bigr)\,\rmd t
\\
&=&\int_{t\ge0} \bigl(E\bigl(e^{-aL^{-1}_1+\ga H_1};1<L_\infty
\bigr) \bigr)^t \,\rmd t<\infty,
\end{eqnarray*}
showing that $\eta_\alpha^V$ is finite on compact sets.

%
\begin{thmm}\label{JOS1} The limiting distribution in \eqref{ltpto}
may be written alternatively as
\begin{eqnarray}
\label{ltpto1} && P\bigl(W_{\tau_0(W)}\in\rmd
\gamma, \tau+\tau_0(W)\in\rmd t\bigr)
\nonumber\\
&&\qquad=\frac{\alpha}{B(T)} \biggl(\mu(\rmd t)e^{-\alpha\gamma}\,\rmd\gamma+
\rmd_{ H}\bigl(\mu*\eta_\alpha^V\bigr) (\rmd t)
\delta_{\{0\}}(\rmd\gamma )\\
&&\hspace*{42pt}\qquad{}+
 \int_{0\le s\le t}\bigl(\mu*\eta_\alpha^V
\bigr) (\rmd t-s) \int_{y\ge0} e^{\alpha y}
\Pi_{{ L}^{-1}, { H}}(\rmd s,y+\rmd \gamma)\,\rmd y \biggr)\nonumber
\end{eqnarray}
for $\gamma\ge0$ and $0\le t< T$.
\end{thmm}

\begin{pf}
By~\cite{GM1}, Corollary 3.1, for $z>0$ and $\gamma, r \ge0$,
%
\begin{eqnarray}
\label{pto1} && P_0\bigl(X_{\tau(z)}-z
\in\rmd\gamma; \tau(z)\in\rmd r\bigr)
\nonumber\\
&&\qquad=\rmd_{ H}\frac{\partial_-}{\partial_- z}V(\rmd r,z)\delta_{\{0\}
}(\rmd
\gamma)
\\
&&\qquad\quad{} +I(\gamma>0)\int_{0\le s\le r}\int_{0\le y\le z}V(\rmd
s,z-\rmd y)\Pi_{{ L}^{-1}, { H}}(\rmd r-s,y+\rmd\gamma),\nonumber
\end{eqnarray}
where ${\partial_-}/{\partial_- z}$ denotes the left derivative. A
straightforward calculation involving changes of variable and orders of
integration, shows that
\begin{eqnarray}
\label{pto2} && I(\gamma>0)\int_{z\ge0}
\alpha e^{\alpha z}\,\rmd z \int_{0\le y\le
z}V(\rmd s,z-\rmd y)
\Pi_{{ L}^{-1}, { H}}(\rmd r-s,y+\rmd\gamma)
\nonumber\\
&&\qquad= I(\gamma>0)\alpha\int_{\zeta\ge0}e^{\alpha\zeta}V(\rmd s,\rmd
\zeta) \int_{y\ge0} e^{\alpha y}\Pi_{{ L}^{-1}, { H}}(\rmd
r-s,y+\rmd\gamma)\,\rmd y
\\
&&\qquad= \alpha\eta_\alpha^V(\rmd s) \int_{y\ge0}
e^{\alpha y}\Pi_{{
L}^{-1}, { H}}(\rmd r-s,y+\rmd\gamma)\,\rmd y,\nonumber
\end{eqnarray}
where the $I(\gamma>0)$ term may be omitted in the final expression
since the measure there assigns no mass to the set $\{\gamma=0\}$.
On the other hand, if $X$ creeps, that is $\rmd_{ H}>0$, then from
\cite{GM1}, Theorem 3.1(ii),
\begin{equation}
\label{pto3} \int_{z\ge0}\alpha e^{\alpha z}\,\rmd z
\frac{\partial_-}{\partial_-
z}V(\rmd r,z) =\alpha\eta_\alpha^V(\rmd r).
\end{equation}
Thus substituting \eqref{pto1}, \eqref{pto2} and \eqref{pto3} into
\eqref{ltpto} gives
\begin{eqnarray*} &&\hspace*{-5pt} P\bigl(W_{\tau_0(W)}\in\rmd\gamma, \tau+
\tau_0(W)\in\rmd t\bigr)
\\
&&\hspace*{17pt}=\frac{\alpha}{B(T)} \biggl(\mu(\rmd t)e^{-\alpha\gamma}\,\rmd\gamma+
\rmd_{ H}\bigl(\mu*\eta_\alpha^V\bigr) (\rmd t)
\delta_{\{0\}}(\rmd\gamma)
\\
&&\hspace*{25pt}\qquad\quad{} + \int_{0\le r\le t}\mu(\rmd t-r)\int_{0\le s\le r}
\eta_\alpha^V(\rmd s) \\
&&\hspace*{175pt}{}\times \int_{y\ge0}
e^{\alpha y}\Pi_{{ L}^{-1}, {
H}}(\rmd r-s,y+\rmd\gamma)\,\rmd y \biggr),
\end{eqnarray*}
which is the same as \eqref{ltpto1}.
\end{pf}

Each of the three terms in \eqref{ltpto1} has a clear meaning.
In order to exit by time $T$, the process must take a large jump from a
neighborhood of the origin to a neighborhood of the boundary. The first
term is a consequence of this jump overshooting the boundary. If the
jump undershoots the boundary, then the process crosses the boundary
either by creeping, which leads to the second term, or by taking a
further (small) jump which results in the final term. From this
description we can read off, for example, that the limiting
(sub-)distribution of the time at which the conditioned process creeps
over the boundary is given by $B(T)^{-1}\ga\rmd_{ H}(\mu*\eta_\alpha^V)$.

The marginal distributions can be obtained from either \eqref{ltpto}
or \eqref{ltpto1}.
We will focus on the latter, but mention in passing that the expression
for the marginal distribution in $t$ obtained from \eqref{ltpto}, is
actually a simple consequence of Theorem~\ref{THMall}; for $0\le t<T$
\begin{equation}
\label{ruinlt} P\bigl(\tau+\tau_0(W)\le t\bigr)=
\frac{B(t)}{B(T)}=\frac{1}{B(T)}\int_{[0,t)}e^{\psi(\ga)(t-s)}Ee^{\alpha\Xbar_{s}}
\,\rmd s.
\end{equation}
%
By integrating out $\gamma$ in \eqref{ltpto1}, and noting that
\begin{equation}
\label{equdP} \ga\int_{y\ge0} \int_{\gamma\ge0}e^{\alpha y}
\Pi_{{ L}^{-1}, {
H}}(\rmd s,y+\rmd\gamma)\,\rmd y=\eta_\ga^{ \Pi_{{ L}^{-1}, {
H}}}(
\rmd s),
\end{equation}
we obtain the alternative characterization
\begin{eqnarray}
\label{ruinf}&& P\bigl(\tau+\tau_0(W)\in\rmd t\bigr)
\nonumber
\\[-8pt]
\\[-8pt]
\nonumber
&&\qquad=
\frac{1}{B(T)} \bigl(\mu(\rmd t) + \ga\rmd_{ H}\bigl(\mu*
\eta_\alpha^V\bigr) (\rmd t)+\bigl(\mu*\eta_\alpha^V*
\eta_\ga^{
\Pi_{{ L}^{-1}, { H}}}\bigr) (\rmd t) \bigr)
\end{eqnarray}
for $0\le t<T$.
Similarly, the marginal distribution in $\gamma$ obtained from~\eqref{ltpto1}~is
%
\begin{eqnarray}\qquad
\label{ltov} &&P(W_{\tau_0(W)}\in\rmd\gamma)\nonumber\\
&&\qquad=
\frac{\alpha}{B(T)} \biggl(\mu (T)e^{-\alpha\gamma}\,\rmd\gamma+
\rmd_{ H}\bigl(\mu*\eta_\alpha^V\bigr) (T)
\delta_{\{0\}}(\rmd\gamma)
\\
&&\qquad\quad\hspace*{32pt}{} + \int_{0\le s<T}\bigl(\mu*\eta_\alpha^V
\bigr) (T-s) \int_{y\ge0} e^{\alpha y}\Pi_{{ L}^{-1}, { H}}(
\rmd s,y+\rmd\gamma)\,\rmd y \biggr),\nonumber 
\end{eqnarray}
where, recall, for any measure $\eta$ and any $t$, $\eta(t)=\eta([0,t))$.\vadjust{\goodbreak}

When $Ee^{\ga X_1}<1$ the limiting distribution of the overshoot was
found in~\cite{kkm} for the case $T=\infty$, that is, conditional on
$\tau(u)<\infty$. To relate the result in~\cite{kkm} to \eqref
{ltov}, we investigate the limit as $T\to\infty$ in \eqref{ltov}.
For this, we recall \eqref{PK2}, which may be rewritten
\begin{equation}
\label{PK1} Ee^{\ga\Xbar_\infty}=q \eta_\alpha^V(\infty).
\end{equation}

%
\begin{thmm} If $\pibar_X^+\in\mathcal{S}^{(\alpha)}$ and $Ee^{\ga
X_1}<1$, then as $T\to\infty$,
\begin{eqnarray}
\label{ltW}&& P(W_{\tau_0(W)}\in\rmd\gamma)
\nonumber
\\[-8pt]
\\[-8pt]
\nonumber
&&\qquad\to\frac{\ga e^{-\alpha\gamma}\,\rmd
\gamma}{Ee^{\ga\Xbar_\infty}} +
\frac{\ga}{q} \biggl(\rmd_H\delta_{0}(\rmd\gamma)+
\int_{y\ge0} e^{\alpha y} \Pi_{ H}(y+\rmd\gamma)
\,\rmd y \biggr),
\end{eqnarray}
where convergence is in the total variation norm.
\end{thmm}

\begin{pf}
Since $\psi(\ga)<0$ when $Ee^{\ga X_1}<1$, on letting $T\to\infty$
we obtain
\begin{equation}
\label{ltin1} \mu(T)\to\mu(\infty)<\infty,
\end{equation}
while by Lemma~\ref{EXbar}, \eqref{PK1} and monotone convergence
\begin{equation}
\label{ltin1a} \bigl(\mu*\eta_\alpha^V\bigr) (T)=\int
_{0\le t<T} \eta_\alpha^V(T-t)\mu (\rmd t)\to
\eta_\alpha^V(\infty)\mu(\infty)<\infty
\end{equation}
and
\begin{eqnarray}
\label{ltin2} &&\int_{0\le s<T}\bigl(
\mu*\eta_\alpha^V\bigr)(T-s) \int_{y\ge0}
e^{\alpha
y}\Pi_{{ L}^{-1}, { H}}(\rmd s,y+\rmd\gamma)\,\rmd y
\nonumber\\
&&\qquad \uparrow\eta_\alpha^V(\infty)\mu(\infty)\int
_{s\ge0}\int_{y\ge0} e^{\alpha y}
\Pi_{{ L}^{-1}, { H}}(\rmd s,y+\rmd\gamma)\,\rmd y
\\
&&\qquad = \eta_\alpha^V(\infty)\mu(\infty)\int
_{y\ge0} e^{\alpha
y}\Pi_{{ H}}(y+\rmd\gamma)
\rmd y, \nonumber\end{eqnarray}
where convergence is in total variation by monotonicity.
Also, again by monotone convergence
\begin{equation}
\label{ltin3} B(T)=\int_{0\le t<T} Ee^{\ga\Xbar_{T-t}}\mu(\rmd t)
\to Ee^{\ga\Xbar_\infty}\mu(\infty)<\infty. 
\end{equation}
Thus \eqref{ltW} follows by letting $T\to\infty$ in
\eqref{ltov} and using \eqref{PK1}.
\end{pf}

The limiting distribution in \eqref{ltW} agrees with the limiting
distribution of the overshoot conditional on $\tu<\infty$ which was
found in~\cite{kkm} (see also~\cite{DK} and~\cite{GM2}, (7.5)).
Since the only possible atoms in the limiting distributions are at $0$,
it thus follows that
%
\begin{equation}
\label{intlts} \lim_{T\to\infty}\lim_{u\to\infty}\PuT(O_u<x)=
\lim_{u\to\infty
}\lim_{T\to\infty}\PuT(O_u<x)
\end{equation}
for every $x\ge0$, when $Ee^{\ga X_1}<1$.\vadjust{\goodbreak}

It is interesting to note that if $Ee^{\ga X_1}= 1$, then $Ee^{\ga
\Xbar_\infty}=\infty$ [see \eqref{Crinf} below], and hence,
formally the limit in \eqref{ltW} becomes
\begin{equation}
\label{ltWCr} P(W_{\tau_0(W)}\in\rmd\gamma)\to\frac{\ga}{q} \biggl(
\rmd_H\delta_{0}(\rmd\gamma)+ \int_{y\ge0}
e^{\alpha y} \Pi_{ H}(y+\rmd\gamma )\,\rmd y \biggr).
\end{equation}
Under the (minor) additional Cram\'er--Lundberg assumption \eqref
{Cr1}, this again agrees with the limiting distribution of the
overshoot conditional on $\tu<\infty$ (see, e.g.,~\cite{GMCr}) and
so \eqref{intlts} continues to hold. However, the argument given above
is not rigorous in this case as all the limiting quantities in \eqref
{ltin1}--\eqref{ltin3} are infinite and hence, cannot be canceled.
To prove \eqref{ltWCr}, more care needs to be taken with the limiting
operations. To this end, we begin by recalling that $Ee^{\ga X_1}= 1$
is equivalent to $\psi(\ga)=0$. Consequently, by the Wiener--Hopf
factorization,
\[
k\gk(-\ga)\whk(\ga)=-\psi(\ga)=0.
\]
Here $ \whk$ is the Laplace exponent of the descending ladder process
$\whH\ge0$ and $k>0$ is some constant depending on the normalization
of the local times (see, e.g.,~\cite{kypbook}, Theorem 6.16(iv)).
Since $X_t\to-\infty$ a.s. when $Ee^{\ga X_1}=1$, it follows that
$\whH$ is a proper (not killed) subordinator, and hence, $Ee^{-\ga
\whH_1}= e^{-\whk(\ga)}<1$. Thus, if $\psi(\ga)=0$, then $\kappa
(-\ga)=0$, and so by \eqref{BivV} and \eqref{PK1},
\begin{equation}
\label{Crinf} Ee^{\ga\Xbar_\infty}=q\int_x e^{\ga x}
\int_{t=0}^\infty P(H_t\in \rmd x,
t<L_\infty)\,\rmd t=\int_0^\infty
e^{-\kappa(-\ga) t}\,\rmd t=\infty.\hspace*{-35pt}
\end{equation}

%
\begin{lemma} Assume $Ee^{\ga X_1}=1$ and \eqref{Cr1},
then
\begin{equation}
\label{ratB} \lim_{T\to\infty}\frac{B(T-t)}{B(T)}\to1.
\end{equation}
\end{lemma}

\begin{pf} First observe that from \eqref{B}
\begin{equation}
\label{Bmu} \frac{B(T)}{T}=\frac1{T} \int
_{0\le s<T} Ee^{\ga\Xbar_s}\,\rmd s\to \infty
\end{equation}
since $Ee^{\ga\Xbar_\infty}=\infty$. Thus by \eqref{Doob}, for
fixed $t$,
\begin{eqnarray*} \frac1{B(T)}\int_{T-t\le s<T}
Ee^{\ga\Xbar_s}\,\rmd s &\le&\frac
{1}{B(T)}\int_{T-t\le s<T}C(1+s)
\,\rmd s
\\
&\le&\frac{C(1+T)t}{B(T)}\to0
\end{eqnarray*}
as $T\to\infty$. Hence, \eqref{ratB} holds.
\end{pf}

%
\begin{lemma} Assume $\pibar_X^+\in\mathcal{S}^{(\alpha)}$,
$Ee^{\ga X_1}=1$ and \eqref{Cr1}, then
\begin{equation}
\label{Blt} \lim_{T\to\infty}\frac{\mu(T)}{B(T)}=0 \quad\mbox{and}\quad
\lim_{T\to\infty}\frac{(\mu*\eta_\alpha^V)(T)}{B(T)}=\frac1q.\vadjust{\goodbreak}
\end{equation}
\end{lemma}

\begin{pf} The first limit follows immediately from \eqref{Bmu} since
$\mu(T)=T$.
For the second limit, first observe that by \eqref{kappa} and monotone
convergence, as $T\uparrow\infty$,
\begin{equation}
\label{pilt} \qquad\quad\ga\rmd_{ H} +\eta_\ga^{ \Pi_{{ L}^{-1}, { H}}}(T)
\uparrow\ga \rmd_{ H} +\int_{y\ge0}
\bigl(e^{\ga y}-1\bigr)\Pi_{{ H}}(\rmd y)=q-\gk(-\ga )=q
\end{equation}
Now by \eqref{ruinf}
\begin{eqnarray*}B_{T}&=&\mu(T) + \ga\rmd_{ H}
\bigl(\mu*\eta_\alpha^V\bigr) (T)+\bigl(\mu*
\eta_\alpha^V*\eta_\ga^{ \Pi_{{ L}^{-1}, { H}}}\bigr) (T)
\\[-2pt]
&=& \mu(T) + \int_{t<T}\bigl[\ga\rmd_{ H}+
\eta_\ga^{ \Pi_{{ L}^{-1}, {
H}}}(T-t)\bigr]\bigl(\mu*\eta_\alpha^V
\bigr) (\rmd t)
\\[-2pt]
&\le&\mu(T) + \bigl[\ga\rmd_{ H}+\eta_\ga^{ \Pi_{{ L}^{-1}, {
H}}}(T)
\bigr]\bigl(\mu*\eta_\alpha^V\bigr) (T).
\end{eqnarray*}
Thus by \eqref{pilt}
\[
\liminf_{T\to\infty}\frac{(\mu*\eta_\alpha^V)(T)}{B(T)}\ge\frac1{q}.
\]
For the upper bound,
fix $T_0>0$. Then by \eqref{ruinf}
\begin{eqnarray*} B_{T+T_0} 
&\ge&\mu(T+T_0) +
\int_{t<T}\bigl[\ga\rmd_{ H}+
\eta_\ga^{ \Pi_{{
L}^{-1}, { H}}}(T+T_0-t)\bigr]\bigl(\mu*
\eta_\alpha^V\bigr) (\rmd t)
\\[-2pt]
&\ge&\bigl[\ga\rmd_{ H}+\eta_\ga^{ \Pi_{{ L}^{-1}, { H}}}(T_0)
\bigr]\bigl(\mu *\eta_\alpha^V\bigr) (T).
\end{eqnarray*}
Dividing by $B(T)$ and using \eqref{ratB}, gives
\[
\limsup_{T\to\infty}\frac{(\mu*\eta_\alpha^V)(T)}{B(T)}\le\frac 1{\ga\rmd_{ H}+
\eta_\ga^{ \Pi_{{ L}^{-1}, { H}}}(T_0)}.
\]
Now let $T_0\to\infty$.
\end{pf}

%
\begin{thmm} Assume $\pibar_X^+\in\mathcal{S}^{(\alpha)}$, $Ee^{\ga
X_1}=1$ and \eqref{Cr1}, then as \mbox{$T\to\infty$},
\begin{equation}
\label{ltW1} P(W_{\tau_0(W)}\in\rmd\gamma)\to\frac{\ga}{q} \biggl(
\rmd_H\delta_{0}(\rmd\gamma)+ \int_{y\ge0}
e^{\alpha y} \Pi_{ H}(y+\rmd\gamma )\,\rmd y \biggr),
\end{equation}
where convergence is in the total variation norm.
\end{thmm}

\begin{pf} Let
\[
f_T(s)=\frac{I(0\le s<T)(\mu*\eta_\alpha^V)(T-s)}{B(T)}.
\]
Then by \eqref{ratB} and \eqref{Blt}, for fixed $s$,
\begin{equation}
\label{ltf} f_T(s)\to\frac{1}{q},
\end{equation}
while
\begin{equation}
\label{ltf1} \sup_s f_T(s)\le\frac{(\mu*\eta_\alpha^V)(T)}{B(T)}\to
\frac{1}{q}.\vadjust{\goodbreak}
\end{equation}
Now from \eqref{equdP} and \eqref{ltov}, for any Borel set $C\subset
[0,\infty)$,
\begin{eqnarray*}&& \biggl\llvert P(W_{\tau_0(W)}\in C)-
\frac{\ga}{q} \biggl(\rmd_H I_{C}(0)+ \int
_{y\ge0} e^{\alpha y} \Pi_{ H}(y+C)\,\rmd y
\biggr)\biggr\rrvert
\\[-2pt]
&&\qquad \le\frac{\mu(T)}{B(T)} + {\alpha\rmd_{ H}}\biggl\llvert
\frac
{(\mu*\eta_\alpha^V)(T)}{B(T)}-q^{-1}\biggr\rrvert
\\[-2pt]
&&\qquad\quad{} +\int_{s\ge0} \bigl\llvert f_T(s)-q^{-1}
\bigr\rrvert\int_{y\ge
0} \ga e^{\alpha y}
\Pi_{{ L}^{-1}, { H}}(\rmd s, y+C)\,\rmd y
\\[-2pt]
&&\qquad \le\frac{\mu(T)}{B(T)} + {\alpha\rmd_{ H}}\biggl\llvert
\frac
{(\mu*\eta_\alpha^V)(T)}{B(T)}-q^{-1}\biggr\rrvert+\int_{s\ge0}
\bigl|f_T(s)-q^{-1}\bigr|\eta_\alpha^{\Pi_{{ L}^{-1}, { H}}}(\rmd
s).
\end{eqnarray*}
Since $\eta_\alpha^{\Pi_{{ L}^{-1}, { H}}}(\rmd s)$ is a finite
measure, the result follows by taking the supremum over all $C$
and using \eqref{Blt}, \eqref{ltf}, \eqref{ltf1} and bounded convergence.
\end{pf}

When $Ee^{\ga X_1}>1$ it is possible that $q>0$ or $q=0$. In either
case it seems more difficult to obtain an analogue of \eqref{ltW1}, in
part because \eqref{ratB} no longer holds. One case in which the limit
in \eqref{ltW1} can be found is when $X$ is a subordinator, and so
$q=0$. In this case we may take $H=X$, and similar calculations to
those above lead to
\[
P(W_{\tau_0(W)}\in\rmd\gamma)\to\frac{\ga}{\psi(\ga)} \biggl(
\rmd_H\delta_{0}(\rmd\gamma)+ \int_{y\ge0}
e^{\alpha y} \Pi_{
H}(y+\rmd\gamma)\,\rmd y \biggr).
\]
On the other hand, since $X_t\to\infty$ a.s.,
\[
\lim_{T\to\infty}\PuT(O_u\in\rmd\gamma)=P(O_u\in\rmd
\gamma),
\]
and
by standard renewal theory (see, e.g.,~\cite{kypbook}, Theorem 5.7)
\[
\lim_{u\to\infty}P(O_u\in\rmd\gamma)=
\frac{1}{m} \biggl(\rmd_H\delta_{0}(\rmd\gamma)+
\int_{y\ge0} \Pi_{ H}(y+\rmd\gamma)\,\rmd y \biggr),
\]
where $m=EH_1=EX_1$. In particular, this shows that \eqref{intlts} no
longer holds when $Ee^{\ga X_1}>1$.

A similar discussion applies to the first passage time; if $Ee^{\ga
X_1}<1$ then by~\eqref{rt7}, \eqref{tauTu} and \eqref{tauuT}
\begin{equation}
\label{intltsr} \lim_{T\to\infty}\lim_{u\to\infty}\PuT\bigl(\tu<t\bigr)=
\lim_{u\to
\infty}\lim_{T\to\infty}\PuT\bigl(\tu<t\bigr)
\end{equation}
for all $t\ge0$. When $Ee^{\ga X_1}\ge1$, letting $T\to\infty$ in
\eqref{ruinlt} shows that for all $t\ge0$
\[
\lim_{T\to\infty}\lim_{u\to\infty}\PuT\bigl(\tu< t\bigr)=0,
\]
while by Theorem~\ref{THMall} and Remark~\ref{rmkinf}
\[
\lim_{u\to\infty}\lim_{T\to\infty}\PuT\bigl(\tu<t\bigr)=\lim_{u\to
\infty}P\bigl(\tu<t|
\tu<\infty\bigr)=0.
\]
Hence, \eqref{intltsr} is also valid in this case, but in a degenerate sense.

From the calculations presented in this section, it is hopefully
clear that the asymptotic behavior of many other functionals of the
path can be investigated in a similar manner.




\begin{appendix}
\section*{Appendix}\label{sA}

The \hyperref[sA]{Appendix} gives more details on several formulas involving $\mu_K$
and $\nu_x$ where $K,x\in(-\infty,\infty]$. They are first defined
for product sets
$(A\times B\times C)\in\mathcal{F}\otimes\mathcal{B}([0,\infty
))\otimes\mathcal{B}$ by
\[
\mu_K(A\times B\times C) = \int_{t\in B}E
\bigl(e^{\alpha X_{{t}-}}; \Xt\in A, X_{{t}-}\in C, X_{{t}-}<K
\bigr) \,\rmd t
\]
and
\[
\nu_x(A\times B\times C) = \int_{\{z>-x\}\cap C}\alpha
e^{-\alpha z}P_z \bigl(X\in A,\tau(0) \in B \bigr)\,\rmd z
\]
and then extended to measures on $\mathcal{F}\otimes\mathcal
{B}([0,\infty))\otimes\mathcal{B}$.

%
\setcounter{lemma}{0}
\begin{lemma}\label{AL2} For each $T\ge0$ the following equality of
measures holds:
\setcounter{equation}{0}
\begin{equation}
\label{mKeq} I(t<T)\mu_K(\rmd w, \rmd t, \rmd\phi)=I\bigl(
\tau_\Delta(w)<T\bigr)\mu_K(\rmd w, \rmd t, \rmd\phi)
\end{equation}
and
\begin{equation}
\label{nKeq} I(r<T)\nu_x\bigl(\rmd w', \rmd r, \rmd z
\bigr)=I\bigl(\tau_0\bigl(w'\bigr)<T\bigr)
\nu_x\bigl(\rmd w', \rmd r, \rmd z\bigr).
\end{equation}
\end{lemma}

\begin{pf} For \eqref{mKeq}, it suffices to show
\[
\int_{A\times B\times C} I(t<T)\mu_K(\rmd w, \rmd t, \rmd
\phi)=\int_{A\times B\times C}I\bigl(\tau_\Delta(w)<T\bigr)
\mu_K(\rmd w, \rmd t, \rmd \phi)
\]
for every $(A\times B\times C)\in\mathcal{F}\otimes\mathcal
{B}([0,\infty))\otimes\mathcal{B}$. Let $\{w\dvtx \tau_\Delta(w)<T\}
=A_1\in\mathcal{F}$.
Observe that
\[
\Xt\in AA_1 \quad\mbox{iff } \Xt\in A \mbox{ and } t<T.
\]
Thus
\begin{eqnarray*}
&&\int_{A\times B\times C} I\bigl(
\tau_\Delta(w)<T\bigr)\mu_K(\rmd w, \rmd t, \rmd\phi)
\\
&&\qquad=\int_{t\in B}E\bigl(e^{\alpha X_{{t}-}}; \Xt\in
AA_1, X_{{t}-}\in C, X_{{t}-}<K\bigr) \,\rmd t
\\
&&\qquad=\int_{t\in B}I(t<T)E\bigl(e^{\alpha X_{{t}-}}; \Xt\in A,
X_{{t}-}\in C, X_{{t}-}<K\bigr) \,\rmd t
\\
&&\qquad= \mu_K\bigl(A\times\bigl(B\cap[0,T)\bigr)\times C\bigr)
\\
&&\qquad=\int_{A\times B\times C} I(t<T)\mu_K(\rmd w, \rmd t, \rmd
\phi). \end{eqnarray*}
The proof of \eqref{nKeq} is analogous.
\end{pf}

%
\begin{lemma}\label{ApL1} For any nonnegative measurable function
$H\dvtx D\to\R$
\begin{equation}
\label{al1} \int_{D} H(w)\mu_K(\rmd w)= \int
_0^\infty E\bigl(e^{\alpha X_{t-}}H(\Xt
);X_{t-}<K\bigr) \,\rmd t
\end{equation}
and
\begin{equation}
\label{al2} \int_{D} H\bigl(w'\bigr)
\nu_x\bigl(\rmd w'\bigr) =\int_{z>-x}
\ga e^{-\ga z} E_z\bigl(H(X);\tau(0)<\infty\bigr) \,\rmd z.
\end{equation}
\end{lemma}

\begin{pf}
Let $H=1_A$ for $A\in\mathcal{F}$. Then
\begin{eqnarray*}&& \int_0^\infty E
\bigl(e^{\alpha X_{t-}}H(\Xt);X_{t-}<K\bigr)\,\rmd t
\\
&&\qquad=\int_0^\infty E\bigl(e^{\alpha X_{t-}}; \Xt\in
A;,X_{t-}<K\bigr)\,\rmd t
\\
&&\qquad=\mu_K(A)=\int_{D} H(w)\mu_K(
\rmd w). \end{eqnarray*}
Formula \eqref{al1} then follows by standard arguments. The proof of
\eqref{al2} is similar.
\end{pf}

Applying \eqref{al1} to the function $H(w)f(\tau_{\Delta
}(w))g(w_{\tau_{\Delta}-},\wbar_{\tau_{\Delta}-})I(\tau_{\Delta
}(w)<T)$ where $f\dvtx \R\to\R$ and $g\dvtx \R^2\to\R$
are nonnegative Borel functions, gives
\begin{eqnarray}
\label{Am1} &&\int_{D}H(w)f
\bigl(\tau_\Delta(w)\bigr)g(w_{\tau_{\Delta}-},\wbar_{\tau
_{\Delta}-})
\mu^T_K(\rmd w)
\nonumber
\\[-8pt]
\\[-8pt]
\nonumber
&&\qquad=\int_{[0,T)} f(t) E\bigl(e^{\alpha X_{t-}}H(
\Xt)g(X_{t-},\Xbar_{t-}); X_{t-}<K\bigr)\,\rmd t.
\end{eqnarray}
As a special case we obtain
\begin{eqnarray}
\label{Am1a} &&\int_{D}H(w) f
\bigl(\tau_\Delta(w)\bigr)\mu^T_K(\rmd w)
\nonumber\\
&&\qquad=\int_{[0,T)} f(t) E\bigl(e^{\alpha X_{t-}}H(\Xt);
X_{t-}<K\bigr)\,\rmd t
\nonumber
\\[-8pt]
\\[-8pt]
\nonumber
&&\qquad=\int_{[0,T)} f(t) \int_{D} H(w) E
\bigl(e^{\alpha X_{t-}}; \Xt\in\rmd w, X_{t-}<K\bigr)\,\rmd t
\\
&&\qquad=\int_{{D}\otimes[0,\infty)}H(w)f(t) \mu^T_K(\rmd
w, \rmd t).\nonumber
\end{eqnarray}
Similarly,
\begin{eqnarray}
\label{An1} &&\int_{D}H
\bigl(w'\bigr)f\bigl(\tau_0\bigl(w'
\bigr)\bigr)\nu^T_x(\rmd w)
\nonumber\\
&&\qquad=\int_{[0,T)} f(r) \int_{z>-x}\ga
e^{-\ga z} E_z\bigl(H(X);\tau(0)\in \rmd r\bigr)
\\
&&\qquad=\int_{{D}\otimes[0,\infty)}H\bigl(w'\bigr)f(r)
\nu^T_x\bigl(\rmd w', \rmd r\bigr).\nonumber
\end{eqnarray}

%
\begin{lemma} If $f,g\dvtx \R^2\to\R$ are nonnegative Borel functions then
\begin{eqnarray}
\label{ltb12} &&\int_{{D}\otimes{D}} f\bigl(
\tau_\Delta(w), \tau_0\bigl(w'\bigr)\bigr)
g(w_{\tau
_{\Delta}-},\wbar_{\tau_{\Delta}-}) \mu^T_K(\rmd
w)\otimes\nu^T_x\bigl(\rmd w'\bigr)
\nonumber
\\[-8pt]
\\[-8pt]
\nonumber
&&\qquad =\int_{[0,T)}\int_{[0,T)} f(t,r)E
\bigl(e^{\alpha
X_{t-}}g(X_{t-}, \Xbar_{t-});
X_{t-}<K\bigr)\,\rmd t \nu_x( \rmd r).
 \end{eqnarray}
\end{lemma}

\begin{pf} By \eqref{Am1} and \eqref{An1}, both with $H\equiv1$,
\eqref{ltb12} holds if $f(t,r)=f_1(t)f_2(r)$. The general result then
follows by standard arguments.
\end{pf}

Setting $g\equiv1$ gives
\begin{eqnarray}
\label{ltb12a}&& \int_{{D}\otimes{D}} f\bigl(
\tau_\Delta(w), \tau_0\bigl(w'\bigr)\bigr)
\mu_K^T(\rmd w)\otimes\nu_x^T
\bigl(\rmd w'\bigr)
\nonumber
\\[-8pt]
\\[-8pt]
\nonumber
&&\qquad=\int_0^\infty
\int_0^\infty f(t,r)\mu_K^T(
\rmd t)\nu_x^T( \rmd r).
\end{eqnarray}
Taking $g(y_1,y_2)=1_{(-\infty, K)}(y_2)$ yields
\begin{eqnarray}
\label{ltb11} &&\int_{{D}\otimes{D}} f\bigl(
\tau_\Delta(w), \tau_0\bigl(w'\bigr)\bigr) I(
\wbar_{\tau
_\Delta-}<K) \mu_K^T(\rmd w)\otimes
\nu_x^T\bigl(\rmd w'\bigr)
\nonumber\\
&&\qquad =\int_0^\infty\int_{0}^\infty
f(t,r)I(t<T)E\bigl(e^{\alpha X_{t-}}; \Xbar_{t-}<K\bigr) \,\rmd t
\nu_x^T( \rmd r)
\\
&&\qquad =\int_0^\infty\int_0^\infty
f(t,r)\mubar_K^T(\rmd t)\nu_x^T(
\rmd r), \nonumber
\end{eqnarray}
where
\[
\mubar_K^T(\rmd t)=I(t<T)E\bigl(e^{\alpha X_{t-}};
\Xbar_{t-}<K\bigr) \,\rmd t.
\]

%
\begin{lemma} Let $H\dvtx {D}\otimes{D}\to\R$ be nonnegative and
measurable, and $f\dvtx \R^2\to\R$ be nonnegative and Borel, then
\begin{eqnarray}
\label{Hf} \qquad&&\int_{{D}\otimes{D}} H
\bigl(w,w'\bigr)f\bigl(\tau_{\Delta}(w),\tau_0
\bigl(w'\bigr)\bigr) \mu^T_K(\rmd w)
\otimes\nu^T_x\bigl(\rmd w'\bigr)
\nonumber
\\[-8pt]
\\[-8pt]
\nonumber
&&\qquad =\int_{ ({D}\otimes[0,\infty) )\otimes({D}\otimes
[0,\infty) )} H\bigl(w,w'\bigr)f(t,r)
\mu_K^T(\rmd w,\rmd t)\otimes\nu^{T}_x
\bigl( \rmd w', \rmd r\bigr) .
\end{eqnarray}
\end{lemma}

\begin{pf} If $H(w,w')=H_1(w)H_2(w')$ and $f(t,r)=f_1(t)f_2(r)$, then
\eqref{Hf} holds by \eqref{Am1a} and \eqref{An1}.
The result then follows by standard arguments.
\end{pf}
As a special case we obtain
\begin{eqnarray}
\label{AlH} &&\int_{{D}\otimes{D}} H
\bigl(w,w'\bigr)I\bigl(\tau_{\Delta}(w)+ \tau_0
\bigl(w'\bigr)<T\bigr) \mu^T_K(\rmd w)
\otimes\nu^T_x\bigl(\rmd w'\bigr)
\nonumber
\\[-8pt]
\\[-8pt]
\nonumber
&&\qquad =\int_{{D}\otimes{D}} H\bigl(w,w'\bigr)\int
_{[0,\infty
)} \mu_K^T(\rmd w,\rmd t)
\nu^{T-t}_x\bigl( \rmd w'\bigr). \end{eqnarray}
\end{appendix}

\section*{Acknowledgments}
The author would like to
thank Ross Maller for his hospitality and willing ear while portions of
this work were being carried out. Thanks also to an anonymous referee
for some helpful comments.

%

%


\printaddresses


\begin{thebibliography}{32}

\bibitem{AS}
%
\begin{barticle}[mr]
\bauthor{\bsnm{Albin},~\bfnm{J.~M.~P.}\binits{J.~M.~P.}} \AND
\bauthor{\bsnm{Sund{\'e}n},~\bfnm{Mattias}\binits{M.}}
(\byear{2009}).
\btitle{On the asymptotic behaviour of {L}\'evy processes.~{I}.
{S}ubexponential and exponential processes}.
\bjournal{Stochastic Process. Appl.}
\bvolume{119}
\bpages{281--304}.
\bid{doi={10.1016/j.spa.2008.02.004}, issn={0304-4149}, mr={2485028}}
\bptok{imsref}%
\end{barticle}
%
\endbibitem

\bibitem{Arp}
%
\begin{bbook}[mr]
\bauthor{\bsnm{Asmussen},~\bfnm{S{\o}ren}\binits{S.}} \AND
\bauthor{\bsnm{Albrecher},~\bfnm{Hansj{\"o}rg}\binits{H.}}
(\byear{2010}).
\btitle{Ruin Probabilities},
\bedition{2nd} ed.
\bpublisher{World Scientific}, \baddress{Hackensack, NJ}.
\bid{doi={10.1142/9789814282536}, mr={2766220}}
\bptok{imsref}%
\end{bbook}
%
\endbibitem

\bibitem{berm}
%
\begin{barticle}[mr]
\bauthor{\bsnm{Berman},~\bfnm{Simeon~M.}\binits{S.~M.}}
(\byear{1986}).
\btitle{The supremum of a process with stationary independent and symmetric
increments}.
\bjournal{Stochastic Process. Appl.}
\bvolume{23}
\bpages{281--290}.
\bid{doi={10.1016/0304-4149(86)90041-4}, issn={0304-4149}, mr={0876050}}
\bptok{imsref}%
\end{barticle}
%
\endbibitem

\bibitem{bert}
%
\begin{bbook}[mr]
\bauthor{\bsnm{Bertoin},~\bfnm{Jean}\binits{J.}}
(\byear{1996}).
\btitle{L\'evy Processes}.
\bseries{Cambridge Tracts in Mathematics}
\bvolume{121}.
\bpublisher{Cambridge Univ. Press}, \baddress{Cambridge}.
\bid{mr={1406564}}
\bptok{imsref}%
\end{bbook}
%
\endbibitem


\bibitem{BGT}
%
\begin{bbook}[mr]
\bauthor{\bsnm{Bingham},~\bfnm{N.~H.}\binits{N.~H.}},
\bauthor{\bsnm{Goldie},~\bfnm{C.~M.}\binits{C.~M.}} \AND
\bauthor{\bsnm{Teugels},~\bfnm{J.~L.}\binits{J.~L.}}
(\byear{1987}).
\btitle{Regular Variation}.
\bseries{Encyclopedia of Mathematics and Its Applications}
\bvolume{27}.
\bpublisher{Cambridge Univ. Press}, \baddress{Cambridge}.
\bid{mr={0898871}}
\bptok{imsref}%
\end{bbook}
%
\endbibitem

\bibitem{br}
%
\begin{barticle}[mr]
\bauthor{\bsnm{Braverman},~\bfnm{Michael}\binits{M.}}
(\byear{1997}).
\btitle{Suprema and sojourn times of {L}\'evy processes with exponential
tails}.
\bjournal{Stochastic Process. Appl.}
\bvolume{68}
\bpages{265--283}.
\bid{doi={10.1016/S0304-4149(97)00031-8}, issn={0304-4149}, mr={1454836}}
\bptok{imsref}%
\end{barticle}
%
\endbibitem

\bibitem{BS}
%
\begin{barticle}[mr]
\bauthor{\bsnm{Braverman},~\bfnm{Michael}\binits{M.}} \AND
\bauthor{\bsnm{Samorodnitsky},~\bfnm{Gennady}\binits{G.}}
(\byear{1995}).
\btitle{Functionals of infinitely divisible stochastic processes with
exponential tails}.
\bjournal{Stochastic Process. Appl.}
\bvolume{56}
\bpages{207--231}.
\bid{doi={10.1016/0304-4149(94)00074-4}, issn={0304-4149}, mr={1325220}}
\bptok{imsref}%
\end{barticle}
%
\endbibitem

\bibitem{Ch}
%
\begin{barticle}[auto:STB|2012/07/16|09:13:59]
\bauthor{\bsnm{Chistiakov},~\bfnm{V.~P.}\binits{V.~P.}}
(\byear{1964}).
\btitle{A theorem on sums of independent random variables and its application
to branching random processes}.
\bjournal{Theory Probab. Appl.}
\bvolume{9}
\bpages{640--648}.
\bptok{imsref}%
\end{barticle}
%
\endbibitem

\bibitem{CNW}
%
\begin{barticle}[mr]
\bauthor{\bsnm{Chover},~\bfnm{J.}\binits{J.}},
\bauthor{\bsnm{Ney},~\bfnm{P.}\binits{P.}} \AND
\bauthor{\bsnm{Wainger},~\bfnm{S.}\binits{S.}}
(\byear{1973}).
\btitle{Degeneracy properties of subcritical branching processes}.
\bjournal{Ann. Probab.}
\bvolume{1}
\bpages{663--673}.
\bid{mr={0348852}}
\bptok{imsref}%
\end{barticle}
%
\endbibitem

\bibitem{C}
%
\begin{barticle}[mr]
\bauthor{\bsnm{Cline},~\bfnm{Daren B.~H.}\binits{D.~B.~H.}}
(\byear{1986}).
\btitle{Convolution tails, product tails and domains of attraction}.
\bjournal{Probab. Theory Related Fields}
\bvolume{72}
\bpages{529--557}.
\bid{doi={10.1007/BF00344720}, issn={0178-8051}, mr={0847385}}
\bptok{imsref}%
\end{barticle}
%
\endbibitem

\bibitem{doneystf}
%
\begin{bincollection}[mr]
\bauthor{\bsnm{Doney},~\bfnm{Ronald~A.}\binits{R.~A.}}
(\byear{2005}).
\btitle{Fluctuation theory for L\'evy processes}.
In \bbooktitle{Ecole d'Et\'{e} de Probabilit\'{e}s de Saint-Flour XXXV---2005}.
\bseries{Lecture Notes in Math.}
\bvolume{1897}
\bpublisher{Springer}, \baddress{Berlin}.
\bid{mr={2126962}}
\bptok{imsref}%
\end{bincollection}
%
\endbibitem

\bibitem{DK}
%
\begin{barticle}[mr]
\bauthor{\bsnm{Doney},~\bfnm{R.~A.}\binits{R.~A.}} \AND
\bauthor{\bsnm{Kyprianou},~\bfnm{A.~E.}\binits{A.~E.}}
(\byear{2006}).
\btitle{Overshoots and undershoots of {L}\'evy processes}.
\bjournal{Ann. Appl. Probab.}
\bvolume{16}
\bpages{91--106}.
\bid{doi={10.1214/105051605000000647}, issn={1050-5164}, mr={2209337}}
\bptok{imsref}%
\end{barticle}
%
\endbibitem

\bibitem{Dur}
%
\begin{bbook}[mr]
\bauthor{\bsnm{Durrett},~\bfnm{Rick}\binits{R.}}
(\byear{2010}).
\btitle{Probability: Theory and Examples},
\bedition{4th} ed.
\bpublisher{Cambridge Univ. Press}, \baddress{Cambridge}.
\bid{mr={2722836}}
\bptok{imsref}%
\end{bbook}
%
\endbibitem

\bibitem{EG}
%
\begin{barticle}[mr]
\bauthor{\bsnm{Embrechts},~\bfnm{Paul}\binits{P.}} \AND
\bauthor{\bsnm{Goldie},~\bfnm{Charles~M.}\binits{C.~M.}}
(\byear{1982}).
\btitle{On convolution tails}.
\bjournal{Stochastic Process. Appl.}
\bvolume{13}
\bpages{263--278}.
\bid{doi={10.1016/0304-4149(82)90013-8}, issn={0304-4149}, mr={0671036}}
\bptok{imsref}%
\end{barticle}
%
\endbibitem

\bibitem{GM1}
%
\begin{barticle}[mr]
\bauthor{\bsnm{Griffin},~\bfnm{Philip~S.}\binits{P.~S.}} \AND
\bauthor{\bsnm{Maller},~\bfnm{Ross~A.}\binits{R.~A.}}
(\byear{2011}).
\btitle{The time at which a {L}\'evy process creeps}.
\bjournal{Electron.~J. Probab.}
\bvolume{16}
\bpages{2182--2202}.
\bid{doi={10.1214/EJP.v16-945}, issn={1083-6489}, mr={2861671}}
\bptok{imsref}%
\end{barticle}
%
\endbibitem

\bibitem{GM2}
%
\begin{barticle}[auto:STB|2012/07/16|09:13:59]
\bauthor{\bsnm{Griffin},~\bfnm{P.~S.}\binits{P.~S.}} \AND
\bauthor{\bsnm{Maller},~\bfnm{R.~A.}\binits{R.~A.}}
(\byear{2012}).
\btitle{Path decomposition of ruinous behaviour for a~general
L\'evy insurance risk process}.
\bjournal{Ann. Appl. Probab.}
\bvolume{22}
\bpages{1411--1449}.
\bptok{imsref}%
\end{barticle}
%
\endbibitem



\bibitem{GMCr}
%
\begin{barticle}[auto:STB|2012/07/16|09:13:59]
\bauthor{\bsnm{Griffin},~\bfnm{P.~S.}\binits{P.~S.}},
\bauthor{\bsnm{Maller},~\bfnm{R.~A.}\binits{R.~A.}} \AND
\bauthor{\bparticle{van} \bsnm{Schaik},~\bfnm{K.}\binits{K.}}
(\byear{2012}).
\btitle{Asymptotic distributions of the overshoot and
undershoots for the L\'evy insurance risk process in the Cram\'er and
convolution equivalent cases}.
\bjournal{Insurance: Mathematics and Economics}
\bvolume{51}
\bpages{382--392}.
\bptok{imsref}%
\end{barticle}
%
\endbibitem

\bibitem{ht}
%
\begin{barticle}[mr]
\bauthor{\bsnm{Hao},~\bfnm{Xuemiao}\binits{X.}} \AND
\bauthor{\bsnm{Tang},~\bfnm{Qihe}\binits{Q.}}
(\byear{2009}).
\btitle{Asymptotic ruin probabilities of the {L}\'evy insurance model under
periodic taxation}.
\bjournal{Astin Bull.}
\bvolume{39}
\bpages{479--494}.
\bid{doi={10.2143/AST.39.2.2044644}, issn={0515-0361}, mr={2751836}}
\bptok{imsref}%
\end{barticle}
%
\endbibitem

\bibitem{HP}
%
\begin{bbook}[mr]
\bauthor{\bsnm{Hille},~\bfnm{Einar}\binits{E.}} \AND
\bauthor{\bsnm{Phillips},~\bfnm{Ralph~S.}\binits{R.~S.}}
(\byear{1957}).
\btitle{Functional Analysis and Semi-groups}.
\bseries{American Mathematical Society Colloquium Publications}
\bvolume{31}.
\bpublisher{Amer. Math. Soc.}, \baddress{Providence, RI.}
\bid{mr={0089373}}
\bptnote{check year}%
\bptok{imsref}%
\end{bbook}
%
\endbibitem

\bibitem{HK}
%
\begin{bincollection}[mr]
\bauthor{\bsnm{Hubalek},~\bfnm{F.}\binits{F.}} \AND
\bauthor{\bsnm{Kyprianou},~\bfnm{E.}\binits{E.}}
(\byear{2011}).
\btitle{Old and new examples of scale functions for spectrally negative
{L}\'evy processes}.
In \bbooktitle{Seminar on {S}tochastic {A}nalysis, {R}andom {F}ields and
{A}pplications {VI}}
(\beditor{\binits{R.} \bsnm{Dalang}},
\beditor{\binits{M.} \bsnm{Dozzi}}
\AND
\beditor{\binits{F.} \bsnm{Russo}}, eds.).
\bseries{Progress in Probability}
\bvolume{63}
\bpages{119--145}.
\bpublisher{Birkh\"auser/Springer Basel AG}, \baddress{Basel}.
\bid{doi={10.1007/978-3-0348-0021-1_8}, mr={2857022}}
\bptnote{check year}%
\bptok{imsref}%
\end{bincollection}
%
\endbibitem

\bibitem{kl}
%
\begin{barticle}[mr]
\bauthor{\bsnm{Kl{\"u}ppelberg},~\bfnm{Claudia}\binits{C.}}
(\byear{1989}).
\btitle{Subexponential distributions and characterizations of related classes}.
\bjournal{Probab. Theory Related Fields}
\bvolume{82}
\bpages{259--269}.
\bid{doi={10.1007/BF00354763}, issn={0178-8051}, mr={0998934}}
\bptok{imsref}%
\end{barticle}
%
\endbibitem

\bibitem{kkm}
%
\begin{barticle}[mr]
\bauthor{\bsnm{Kl{\"u}ppelberg},~\bfnm{Claudia}\binits{C.}},
\bauthor{\bsnm{Kyprianou},~\bfnm{Andreas~E.}\binits{A.~E.}} \AND
\bauthor{\bsnm{Maller},~\bfnm{Ross~A.}\binits{R.~A.}}
(\byear{2004}).
\btitle{Ruin probabilities and overshoots for general {L}\'evy
insurance risk
processes}.
\bjournal{Ann. Appl. Probab.}
\bvolume{14}
\bpages{1766--1801}.
\bid{doi={10.1214/105051604000000927}, issn={1050-5164}, mr={2099651}}
\bptok{imsref}%
\end{barticle}
%
\endbibitem

\bibitem{kypbook}
%
\begin{bbook}[mr]
\bauthor{\bsnm{Kyprianou},~\bfnm{Andreas~E.}\binits{A.~E.}}
(\byear{2006}).
\btitle{Introductory Lectures on Fluctuations of {L}\'evy Processes with
Applications}.
\bpublisher{Springer}, \baddress{Berlin}.
\bid{mr={2250061}}
\bptnote{check year}%
\bptok{imsref}%
\end{bbook}
%
\endbibitem

\bibitem{M}
%
\begin{barticle}[mr]
\bauthor{\bsnm{Marcus},~\bfnm{Michael~B.}\binits{M.~B.}}
(\byear{1987}).
\btitle{{$\xi$}-radial processes and random {F}ourier series}.
\bjournal{Mem. Amer. Math. Soc.}
\bvolume{68}
\bpages{viii+181}.
\bid{issn={0065-9266}, mr={0897272}}
\bptok{imsref}%
\end{barticle}
%
\endbibitem

\bibitem{Pakes}
%
\begin{barticle}[mr]
\bauthor{\bsnm{Pakes},~\bfnm{Anthony~G.}\binits{A.~G.}}
(\byear{2004}).
\btitle{Convolution equivalence and infinite divisibility}.
\bjournal{J. Appl. Probab.}
\bvolume{41}
\bpages{407--424}.
\bid{issn={0021-9002}, mr={2052581}}
\bptok{imsref}%
\end{barticle}
%
\endbibitem

\bibitem{Pakes2}
%
\begin{barticle}[mr]
\bauthor{\bsnm{Pakes},~\bfnm{Anthony~G.}\binits{A.~G.}}
(\byear{2007}).
\btitle{Convolution equivalence and infinite divisibility: Corrections and
corollaries}.
\bjournal{J. Appl. Probab.}
\bvolume{44}
\bpages{295--305}.
\bid{doi={10.1239/jap/1183667402}, issn={0021-9002}, mr={2340199}}
\bptok{imsref}%
\end{barticle}
%
\endbibitem

\bibitem{RS}
%
\begin{barticle}[mr]
\bauthor{\bsnm{Rosi{\'n}ski},~\bfnm{Jan}\binits{J.}} \AND
\bauthor{\bsnm{Samorodnitsky},~\bfnm{Gennady}\binits{G.}}
(\byear{1993}).
\btitle{Distributions of subadditive functionals of sample paths of infinitely
divisible processes}.
\bjournal{Ann. Probab.}
\bvolume{21}
\bpages{996--1014}.
\bid{issn={0091-1798}, mr={1217577}}
\bptok{imsref}%
\end{barticle}
%
\endbibitem

\bibitem{sato}
%
\begin{bbook}[mr]
\bauthor{\bsnm{Sato},~\bfnm{Ken-iti}\binits{K.-i.}}
(\byear{1999}).
\btitle{L\'evy Processes and Infinitely Divisible Distributions}.
\bseries{Cambridge Studies in Advanced Mathematics}
\bvolume{68}.
\bpublisher{Cambridge Univ. Press}, \baddress{Cambridge}.
\bid{mr={1739520}}
\bptok{imsref}%
\end{bbook}
%
\endbibitem

\bibitem{Sb}
%
\begin{barticle}[mr]
\bauthor{\bsnm{Sgibnev},~\bfnm{M.~S.}\binits{M.~S.}}
(\byear{1990}).
\btitle{The asymptotics of infinitely divisible distributions in
{$\mathbf{R}$}}.
\bjournal{Sibirsk. Mat. Zh.}
\bvolume{31}
\bpages{135--140, 221}.
\bid{doi={10.1007/BF00971156}, issn={0037-4474}, mr={1046818}}
\bptok{imsref}%
\end{barticle}
%
\endbibitem

\bibitem{Vig}
%
\begin{barticle}[mr]
\bauthor{\bsnm{Vigon},~\bfnm{Vincent}\binits{V.}}
(\byear{2002}).
\btitle{Votre {L}\'evy rampe-t-il?}
\bjournal{J. Lond. Math. Soc. (2)}
\bvolume{65}
\bpages{243--256}.
\bid{doi={10.1112/S0024610701002885}, issn={0024-6107}, mr={1875147}}
\bptok{imsref}%
\end{barticle}
%
\endbibitem

\end{thebibliography}
\end{document}